\documentclass[a4paper,notitlepage, 12pt,reqno]{article}

  \usepackage{mathtext}
 \usepackage[cp1251]{inputenc}
  \usepackage[T2A]{fontenc}
 \usepackage[russian]{babel}
 \usepackage[all,arc,poly,2cell,curve,arrow,tips]{xy}

  \usepackage{enumerate, eucal, amsthm,amsmath, amssymb}

  \allowdisplaybreaks[4]

 \theoremstyle{definition}
 \newtheorem{defn}{Definition}%[section]

 \theoremstyle{plain}
 \newtheorem{thm}{Theorem}
 \newtheorem*{thm*}{Tеорема}
 \newtheorem{prop}{Proposition}
  \newtheorem*{prop*}{Предложение}
 \newtheorem{cor}{Corrolary}
  \newtheorem*{cor*}{Следствие}
 \newtheorem{lem}{Lemma}
  \newtheorem*{lem*}{Лемма}

 \theoremstyle{remark}

 \newtheorem*{remark*}{Замечание}

  \newcounter{ab}
\title{}

\title{ $W$-algebras and higher analogs of the Kniznik-Zamolodchikov equations.\footnote
 {The work was supported by grants NSch-5998.2012.1, RFFI-12-01-31414, MK-4594.2013.1.}}

 \author{D.V. Artamonov\footnote{artamonov.dmitri@gmail.com}, V.A. Golubeva\footnote{goloubeva@yahoo.com}}

\date{}

\begin{document}

\maketitle

\begin{abstract}
The key role in the derivation of the Knizhnik-Zamolodchikov
equations in the $WZW$-theory is played by the energy-momentum
tensor, that is constructed from a central Casimir element of the
second order in a universal enveloping algebra of a corresponding
Lie algebra. In the paper a possibility of construction of analogs
of Knizhnik-Zamolodchikov equations using higher order central
elements is investigated. The  Gelfand elements of the third order
for a simple Lie algebra of series $A$ and Capelli elements of the
fourth order for the a simple Lie algebra of series $B$, $D$ are
considered. In the first case the construction is not possible a the
second case the desired equation is derived.

\end{abstract}

\section{Introduction}

In the $WZW$ theory, associated with the Lie algebra $\mathfrak{g}$
it is proved that the correlation functions of $WZW$-primary fields
satisfy a system of differential equations that is called the
Knizhnik-Zamolodchikov equations.  One obtains this system when one
equates two representation of the action of the Virasoro operator
$L_{-1}$: one as a differential operator and the other as a matrix
operator. The operator $L_{-1}$ is defined by an expansion of the
energy-momentum tensor

\begin{equation}
T(z)=\frac{1}{2(k+g)}\sum_{\alpha}(J^{\alpha}J^{\alpha})(z)=\sum_kL_kz^{-k-1},
\end{equation}

where $J^{\alpha}$  is an orthonomal base of $\mathfrak{g}$ with
respect to the Killing form and $J^{\alpha}(z)$ is a current
corresponding to $J^{\alpha}$. Also $g$ is a dual Coxeter number,
and $k$ is a constant from the $WZW$-action.
 In the derivation of the Knizhnik-Zamolodchikov equation a key role is played
 by the fact that the element $\sum_{\alpha} J^{\alpha}J^{\alpha}\in
U(\mathfrak{g})$ is central.

A natural question arises: is it possible to obtain new equations
using analogs of the energy-momentum tensor that are constructed
from central elements of higher orders?

As a first step in the derivation of an analog of the
Knizhnik-Zamolodchikov equation in the case of a central element $
W=\sum d^{\alpha_1,...,\alpha_n} J^{\alpha_1}...J^{\alpha_n}$ of the
higher order one constructs a field, that is a higher analog of the
energy-mometum tensor

$$W(z)=\sum
d^{\alpha_1,...,\alpha_n} (J^{\alpha_1}(...J^{\alpha_n}))(z)=\sum_k W_kz^{-n-k}.$$

An algebra generates by elements $1,L_n,W_m$ is called the
$W$-algebra \cite{b2}.

As the second step a class of  considered fields is fixed and the
action of
 $W_{-1}$  on these fields is represented as a differential operator. For this an operator product expansion
 of $W(z)$ and the considered field $\varphi(w)$ is investigated. To able to represent the action
  of $W_{-1}$ onto  $\varphi(w)$ as a differential operator this
  expansion must be of type \footnote{In the paper only
singular terms in the OPE are written}

\begin{equation}\label{opeosn}
W(z)\varphi(w)=\frac{const\varphi(w)}{(z-w)^n}+\frac{\mathcal{D}\partial
\varphi(w)}{(z-w)^{n-1}}+...,
\end{equation}
where $\mathcal{D}$ is a differential operator.

In the case of the energy-momentum tensor $T(z)$  of the second
order is postulated by the Ward identity that has geometric origin.

In the present paper instead of $WZW$-primary fields the currents
$J^{\alpha}(z)$  are used. Note that if $k=0$ then currents
$J^{\alpha}(z)$ are $WZW$-primary fields.  Then formulas for the
operator expansions of $T(z)$ and $J^{\alpha}(w)$ is obtained from
the relation of the operator algebra.
%The case of an arbitrary
%tensor representation is considered by taking  a tensor power of an
%adjoint representation and then passing to irreducible components.

This algebraic approach is used for the construction of the OPE of
the energy-momentum tensor constructed using a Casimir element of
the higher order $ W(z)$ and a considered field $\varphi(w)$ which
is constructed from currents.

As the third step  the action of  $W_{-1}$  is represented as an
algebraic operator. But this can be easily done.

Below  we try to realize this construction using  two central
elements of higher orders: the Gelfand element of the third order in
the case when $\mathfrak{g}$ is a simple Lie algebra that belongs to
the series $A$, and the Capelli element of the fourth order  in the
case when $\mathfrak{g}$ belongs to the series $B$, or $D$.

The case  of the energy-momentum tensor associated with the Gelfand
element of the third order
$W=d^{\alpha,\beta,\gamma}J^{\alpha}J^{\beta}J^{\gamma}$ for the
series $A$ is considered in Section \ref{kzs1}. It is shown that it
is not possible to construct the Knizhnik-Zamolodchikov equations.

The following notation is made. The OPE of the current
$J^{\alpha}(z)$ and the field
$W(w)=d^{\alpha,\beta,\gamma}(J^{\alpha}(J^{\beta}J^{\gamma}))(w)$
is of type

\begin{equation}
J^{\alpha}(z)W(w)=\frac{1}{(z-w)^3}d^{\alpha,\beta,\gamma}(J^{\beta}J^{\gamma})(w)\footnote{A
summation over repeating indices is suggested.}.
\end{equation}

Thus it is not an OPE of type \eqref{opeosn}: when one takes the OPE
the field  $J^{\alpha}(z)$ removes  $J^{\alpha}$ from the normal
order product $(J^{\alpha}(J^{\beta}J^{\gamma}))(w)$ and one obtains
$(J^{\beta}J^{\gamma})(w)$.

From this notation the following conjecture arises. To obtain the
OPE of type \eqref{opeosn}, it is necessary  to take a field
$\varphi^{\alpha}(z)$ and a higher analog of the energy momentum
tensor of type
$W(w)=\sum_{\alpha}(\varphi^{\alpha}\varphi^{\alpha})(w)$. In
particular the corresponding central element is a sum of squares
$\sum_{\alpha}\varphi^{\alpha}\varphi^{\alpha}$.

Such central elements there exist in the universal enveloping
algebra of the orthogonal algebra. These are the Capelli elements,
which are sums of squares of noncommutative pfaffians.

This conjecture is verified in Section \ref{kzs2}. As $W$ the
Capelli element of the fourth order is taken, which is a sum of
squares of noncommutative pfaffians and as $\varphi^{\alpha}(z)$
 a noncommutive pfaffian.  It is shown that in these settings
 the construction of an analog of the Knizhnik-Zamolodchikov equation is possible.

\subsection{The content of the paper}

In Section \ref{pred}  we give some preliminary facts about affine
algebras, construction of primary  fields in the  $WZW$-theory,
extensions of the Virasoro algebra, higher order central elements in
the universal enveloping algebra of a simple Lie algebra.

In Section \ref{kzs1} we investigate the  question of the
possibility of the construction of an analog of the
Knizhnik-Zamolodchikov equation using higher order energy-momentum
tensor that is associated with the Gelfand element of the third
element for the series   $A$. The answer is negative.

In Sections \ref{kzs2}, \ref{ur} it is shown that if one takes the
energy-momentum tensor that is associated with the Capelli element
of the fourth order the answer is negative.

Some technical details can be found in  Appendix \ref{ap}.

\section{Preliminaries}
\label{pred}

\subsection{Affine algebras. Currents}
\label{pred1} Let $\mathfrak{g}$ be a simple Lie algebra, denote as
 $\hat{\mathfrak{g}}$ the corresponding non-twisted affine Lie algebra. If $J^{\alpha}$ is a base of $\mathfrak{g}$, then the base of the untwisted
 affine Lie algebra $J^{\alpha}_n$, $K$,
$n\in\mathbb{Z}$.

Define a power series

\begin{equation}
J^{\alpha}(z)=\sum_n J^{\alpha}_n
z^{-n-1}.
\end{equation}
Then one obtains that

\begin{equation}
 \hat{\mathfrak{g}}=(\mathfrak{g}\otimes
\mathbb{C}[[t^{-1},t]])\oplus \mathbb{C}K. \end{equation}

The defining commutations  relation on the language of power series
are written as follows

\begin{align}
\begin{split}
%& \hat{\mathfrak{g}}=(\mathfrak{g}\otimes \mathbb{C}[[t^{-1},t]])\oplus \mathbb{C}K,\\
& [a(t)\oplus nK,b(t)\oplus mK]=[a(t),b(t)]\oplus \omega(a,b)K,
\end{split}
\end{align}
 where $a(t),b(t)\in \mathfrak{g}\otimes \mathbb{C}[[t^{-1},t]]$,
 and
 $\omega$ is the Killing form. In particular the element $K$ is central. The power series $a(t)\in \mathfrak{g}\otimes
 \mathbb{C}[[t^{-1},t]]$ is called a current.

We call a field an arbitrary power series  $U(\mathfrak{g})\otimes
\mathbb{C}[[t^{-1},t]]$, that is a power series with coefficients in
$U(\mathfrak{g})$.

 Write a product of power series  $a(z)$ and $b(w)$ as follows

 \begin{equation}
 a(z)b(w)=\sum_k \frac{(ab)_k(w)}{(z-w)^k}.
 \end{equation}

The coefficient $$(ab)_0(w)=:(ab)(w)$$ is called a normal orderd
product.

\subsection{Currents and $WZW$-primary fields}
\label{prim}

Below the following explicit construction of  $WZW$-primary fields
is used.

Let  $\hat{\mathfrak{g}}=<J_n^a,K>$  be a nontwisted affine algebra
and $\hat{\hat{\mathfrak{g}}}=<\varphi_n^a>$ -- be a loop algebra.
The projection $\pi:\hat{\mathfrak{g}}\rightarrow
\hat{\mathfrak{g}}$, that maps $J_n^a$ to $\varphi_n^a$, and   $K$
into zero, is a homomorphism of Lie algebras.  Thus the algebra
$\hat{\mathfrak{g}}$ act onto $\hat{\hat{\mathfrak{g}}}$ in an
adjoint way with the level  $0$.

Put

$$J^a(z)=\sum_n J_n^a z^{-n-1},\,\,\, \varphi^a(z)=\sum_n\varphi_n^a z^{-n-1}.$$

Then the following OPE takes place

$$J^a(z)\varphi^b(w)=\frac{f_{a,b,c}\varphi_c(w)}{(z-w)}.$$

Thus   $\varphi^b(w)$ are primary fields that transform in the
adjoint representation.

Using fields $\varphi^b(w)$ one easily obtains fields that transform
in the tensor power of the adjoint representation and then as an
irreducible components of this representation.

 In the case of a simple Lie algebra of type   $A$, $B$, $D$ in such a way
 one constructs field that transform in an arbitrary tensor  representation.

%Everywhere below a  $WZW$ theory associated with an adjoint
%representation is considered. Actually we do not put $k=0$, all
%considerations are done for arbitrary  $k$.

Below the fields $J^a(z)$ are considered for arbitrary $k$.

\subsection{The Virasoro algebra and it's extensions. The energy momentum tensor}

In the  $WZW$-theory the following fact takes place that is called
the Sugawara construction \cite{Di}. Suggest that the base
$J^{\alpha}$  is orthonomal with respect to the Killing form.
Consider the energy-momentum tensor

\begin{equation}
T(z)=\frac{1}{2(k+g)}\sum_{\alpha}(J^{\alpha}J^{\alpha})(z),
\end{equation}

and consider the expansion

\begin{equation}
T(z)=\sum_n z^{-n-2}L_n
\end{equation}

then  $L_n$ form the Virasoro algebra.

In the theory  $WZW$  it is natural to consider extensions of the
Virasoro algebra   - the Casimir $W$-algebras. One obtains them when
adds to the energy-momentum tensor one more field that is
constructed from a Casimir element of a higher order. See
description of these $W$-algebras in Sections \ref{wa}, \ref{wb}.

\subsection{Higher order central elements}

 Below two central elements in the universal enveloping algebra
  $U(\mathfrak{g})$ of a simple Lie algebra $\mathfrak{g}$ are defined. They are a Gelfand element of the
  third order in the case of the series  $A$ and the Capelli element in the case of the series $B$ and $D$.

\subsubsection{Gelfand element of the third order for the series $A$.}
\label{s3}

Consider the algebra $\mathfrak{sl}_N$.  Chose as a base the
generalized Gell-Mann matrices $\lambda_{\alpha}$, $\alpha=1,...,N$,
then one has

\begin{equation}
d^{\alpha,\beta,\gamma}=Tr(\lambda_{\alpha}\{\lambda_{\beta},\lambda_{\gamma}\}),\,\,\,\,
\{\lambda_{\beta},\lambda_{\gamma}\}=\lambda_{\beta}\lambda_{\gamma}+\lambda_{\gamma}\lambda_{\beta}.
\end{equation}

Define the third order Gelfand element by the formula

\begin{equation}\label{qq}
W=d^{\alpha,\beta,\gamma}\lambda_{\alpha}\lambda_{\beta}\lambda_{\gamma}.
\end{equation}

%This formulas is used in \cite{b1}. Give another one. Define a
%matrix
%\begin{equation}
%\Lambda=\sum_i \lambda_i\otimes \lambda_i.
%\end{equation}

%\begin{lem} \label{la}
%One has
%$$W=sTr_1(\Lambda^3),$$
%where $sTr_1$  is a symmetric trace in the first factor.
%\end{lem}

%\proof

%The expression  \eqref{qq} for the element $W$ can be written as
%follows

%\begin{equation}
%W=Tr_1(\lambda_{\alpha}\lambda_{\beta}\lambda_{\gamma}\otimes
%\lambda_{\alpha}\lambda_{\beta}\lambda_{\gamma}+\lambda_{\alpha}\lambda_{\beta}\lambda_{\gamma}\otimes
%\lambda_{\alpha}\lambda_{\gamma}\lambda_{\beta}),
%\end{equation}

%where $Tr_1$  is a trace in the first tensor factor. This expression
%equals $sTr_1(\Lambda^3)$.

%\endproof

\subsubsection{Noncommutative pfaffians and Capelli elements}
\label{pred2}

Define central elements in the universal enveloping algebra for
$\mathfrak{o}_N$ that is for a simple Lie algebra of type
 $B$ or $D$.  These elements are constructed using noncommutative pfaffians.  Let us give their
 definition.

Let the algebra $\mathfrak{o}_N$ be realized as the algebra of
skew-symmetric matrices.  In this realization the algebra has a base
that consists of matrices $F_{ij}=E_{ij}-E_{ji}$, $i,j=1,...,n$,
$i<j$. Here $E_{ij}$ are matrix units. This base is orthogonal with
respect to the Killing form.

Commutation relations between these matrices are the following:

\begin{equation}
[F_{ij},F_{kl}]=\delta_{kj}F_{il}-\delta_{il}F_{kj}-\delta_{ik}F_{jl}+\delta_{jl}F_{ki}.
\end{equation}

Let  $\Phi=(\Phi_{ij})$, $i,j=1,...,2k$ be a skew-symmetric
 $2k\times 2k$-matrix,  whose matrix elements belong to a noncommutative ring.

\begin{defn}The noncommutative pfaffian of  $\Phi$ is defined by the
formula

\begin{equation}
Pf \Phi=\frac{1}{k!2^k}\sum_{\sigma\in
S_{2k}}(-1)^{\sigma}\Phi_{\sigma(1)\sigma(2)}...\Phi_{\sigma(2k-1)\sigma(2k)},
\end{equation}

where $\sigma$ is a permutation of the set  $\{1,...,2k\}$.
\end{defn}

Define a matrix $F=(F_{ij})$, $i,j=1,...,N$. Since$F_{ij}=-F_{ji}$,
then $F$ is skew symmetric.

Below we consider pfaffians in $U(\mathfrak{o}_N)$ of the matrix $F$
and matrices $F_I$, that are constructed as follows. Let
$I\subset\{1,...,N\}$ be  a set of indices.  Denote as $F_I$ a
submatrix in $F$ whose rows and columns are indexed by elements from
the set $I$.

The following theorem takes place.

\begin{thm} (\cite{Molev}, $\S$ 7.6)
Let

\begin{equation}
C_k=\sum_{|I|=k,I\subset\{1,...,N\}}(PfF_I)^2,\,\,\,\,k=2,4,...,2[\frac{N}{2}].\end{equation}
 Then $C_k$ belong to the center of  $U(\mathfrak{o}_N)$. In the case of odd $N$ they are algebraically independent and generate the center.
 In the case of even
 $N$  the same is true if one takes instead of $C_N=(PfF)^2$ the element $PfF$.
\end{thm}

The elements $C_k$ are called the Capelli elements.

In the paper \cite{AG} the commutation relations between pfaffians
$PfF_{I}$ and generators were found. To formulate them let us define
$F_{ij}I$. Let  $I=\{i_1,...,i_k\}$, $i_r\in\{1,...,N\}$ be a set of
indices. Identify $i_r$ with the vector $e_{i_r}$ of the standard
representation $V$ of the algebra $\mathfrak{o}_N$, identify the set
$I=\{i_1,...,i_k\}$ with the tensor
 $e_{i_1}\otimes ...\otimes e_{i_k}\in V^{\otimes k}$. Then $F_{ij}I$ is defined as the tensor
 that is obtained when one applies
  $F_{ij}$ to the tensor  $I$.

For numbers $\alpha,\beta\in\mathbb{C}$ define $$PfF_{\alpha I+\beta
J}:=\alpha PfF_I+\beta Pf F_J.$$ Then for every $g\in\mathfrak{o}_N$
the expression $PfF_{gI}$ is well-defined.

\begin{prop}(see \cite{AG})
\label{pffco} $[F_{ij},PfF_I]=PfF_{F_{ij}I}$
\end{prop}

\section{Nonstandard notations}

In the paper some nonstandard notations are widely used. Above one
such notation was introduced
$$PfF_{F_{ij}I}.$$

Let $I$ be divided into subsets $I_1$ and $I_2$. Denote as
$$(-1)^{(I_1,I_2)}$$ the sign of the permutation $I$, which first places  the set  $I_1$, and then the set
 $I_2$.

As above identify the set of indices $I=\{i_1,...,i_k\}$ with a
tensor
$$e_{i_1}\wedge...\wedge e_{i_k}.$$

Define $$I\setminus J$$ as follows. If $J\subset I$, then remove
from $I$ the indices that belong to $J$, construct the corresponding
tensor and take it with the sign $(-1)^{(J,I\setminus J)}$. If
$J\nsubseteq I$, then put $I\setminus J=0$.

Define
$$PfF_{I\setminus J}$$ as a pfaffian of the matrix $(F_{i,j})_{i,j\in I\setminus
J}$ taken with the sign $(-1)^{(J,I\setminus J)}$ in the case
$J\subset I$ and put $PfF_{I\setminus J}=0$ otherwise.

Let as be given a subset $J$ in the set $I$, denote as
$$PfF_{I\setminus J}$$ the pfaffian of the matrix  $(F_{i,j})_{i,j\in I\setminus
J}$, multiplied by the sign $(-1)^{(J,I\setminus J)}$.

For the set $J=\{j_1,j_2\}$ denote as $$F_J$$ the generator
$F_{j_1,j_2}$.

\section{ The Gelfand element of the third order}
\label{kzs1}

Below we investigate the possibility of construction of a
Knizhhnik-Zamolodchikov type equation associated with a Gelfand
central element of the third order.

\subsection{The energy-momentum tensor  of the third order}

In the universal enveloping algebra of a simple Lie  algebra of type
$A_{N+1}$ there exist a Casimir element f the third order

\begin{equation}\label{q}
W= d^{\alpha,\beta,\gamma}J^{\alpha}J^{\beta}J^{\gamma},
\end{equation}

where $d^{\alpha,\beta,\gamma}$ is an invariant traceless tensor.
In \cite{b1}  the following field is defined

\begin{equation}
W(z)=A_N(k)d^{\alpha,\beta,\gamma}(J^{\alpha}(J^{\beta}J^{\gamma}))(z),
\end{equation}

where the normalization constant in the of the algebra  $A_{N+1}$ is
defined by the equality

\begin{equation}
A_N(k)=\sqrt{\frac{N}{18(k+N)^2(N+2k)(N^2-4)}}
\end{equation}

Since the tensor $d^{\alpha,\beta,\gamma}$ is traceless, the element
does not depend  on the placement of brackets in this normal ordered
product.

As the element $T(z)$ the element $W(z)$ is a Sugawara element in
the corresponding untwisted affine Lie algebra   \cite{cht}.

Let us show that there exist no Knizhnik-Zamolodchikov equation
associated with this energy-momentum tensor.

\subsection{The algebra $WA_{2}$ and fields  $W^{\alpha}$}
 \label{wa}

 Consider the $WZW$-theory, associated with a Lie algebra $A_{N+1}$.

 Since the field $W(z)$ has the conformal dimension $3$ it's decomposition is the following

\begin{align}
\begin{split}
W(z)=\sum_{n}(z-w)^{-n-3}W_n(w).
\end{split}
\end{align}

The action of the modes of this decomposition on the field $A(w)$ is
defined by the equality

\begin{align}
\begin{split}\label{f0}
W(z)A(w)=\sum_{n}(z-w)^{-n-3}(W_nA)(w).
\end{split}
\end{align}

\begin{defn}

Introduce a field $W^{\alpha}(z)$

\begin{equation}
W^{\alpha}(z)=\frac{1}{2}d^{\alpha,\beta,\gamma}(J^{\beta}J^{\gamma})(z),
\end{equation}
\end{defn}

In \cite{b1} it is shown that for the fields $W^{\alpha}(z)$ the
following operator expansions take place

\begin{align}
\label{op1}
&T(z)W^a(w)=\frac{1}{(z-w)^2}W^a(w)+\frac{1/2}{z-w}\partial W^a(w),
\end{align}

\begin{align}
\label{ja}
&J^a(z)W^b(w)=\frac{(k+\frac{1}{2}N)}{(z-w)^2}d^{a,b,c}J^c(w)+\frac{1}{z-w}(f_{a,b,c}W^c(w)),
\end{align}

where $f_{a,b,c}$ are structure constants of the algebra.

Thus the fields $W^a(w)$ are Virasoro primary but they are not in
general  $WZW$-primary fields.

\section{Some OPEs}
\label{ope}

In this Section some operator expansions are calculated.

As it is shown in \cite{b1}, \cite{b2}  in addition to operator
expansions \eqref{op1}, \eqref{ja}  the following expansions take
place:

\begin{equation}
\label{r1} J^a(z)W(w)=\frac{k+N}{(z-w)^2}W^a(w).
\end{equation}

Let us find the operator expansion $W(z)W^{\gamma}(w)$.  One has the
following formula for the contraction
$\overbrace{W(z)(J^{\alpha}J^{\beta})(w)}$:

\begin{align}
\begin{split}
& \overbrace{W(z)(J^{\alpha}J^{\beta})(w)}=\frac{1}{2\pi
i}\oint\frac{dx}{x-w} \overbrace{W(z)J^{\alpha}(x)}
J^{\beta}(w)+J^{\alpha}(x)\overbrace{W(z)J^{\beta}(w)}=\\&=\frac{1}{2\pi
i}\oint\frac{dx}{x-w}\frac{(k+N)W^{\alpha}(z)}{(z-x)^2}J^{\beta}(w)+J^{\alpha}(x)\frac{(k+N)W^{\beta}(z)}{(z-w)^2}=\\&=\frac{1}{2\pi
i}\oint\frac{dx}{x-w}\frac{(k+N)}{(z-x)^2}(\frac{(k+\frac{1}{2}N)d^{\alpha,\beta,c}J^c(z)}{(z-w)^2}+\frac{f_{\beta,\alpha,c}W^{c}(z)}{z-w})+\\&+
\frac{(k+N)}{(z-w)^2}(\frac{(k+\frac{1}{2}N)d^{\alpha,\beta,c}J^c(z)}{(x-z)^2}+\frac{f_{\beta,\alpha,c}W^{c}(z)}{x-z})=\\&=
(k+N)(\frac{2(k+\frac{1}{2}N)d^{\alpha,\beta,c}J^c(z)}{(z-w)^4}+\frac{2f_{\beta,\alpha,c}W^c(z)}{(z-w)^3})
\end{split}
\end{align}

From here one obtains

\begin{align}
\begin{split}
\label{r2} &
W(z)W^{\gamma}(w)=\frac{1}{2}d^{\gamma,\alpha,\beta}(k+N)
(\frac{2(k+\frac{1}{2}N)d^{\alpha,\beta,c}J^c(z)}{(z-w)^4}+\frac{f_{\beta,\alpha,c}W^c(z)}{(z-w)^3})=\\&=
(k+N)(k+\frac{N}{2})\frac{2}{N}(N^2-4)\frac{J^{\gamma}(z)}{(z-w)^4}.
%=\\
%&=(k+N)(k+\frac{N}{2})\frac{2}{N}(N^2-4)(\frac{J^{\gamma}(w)}{(z-w)^4}+
%\frac{\partial
%J^{\gamma}(w)}{(z-w)^3}+\frac{\partial^2J^{\gamma}(w)}{(z-w)^2}+\frac{\partial^3
%J^{\gamma}(w)}{(z-w)})
\end{split}
\end{align}

\section{The action of the operators $W_{n}$ on the fields $J^a(z)$, $W^a(z)$ and their derivatives}
\label{dif}

From the formula \eqref{r1} it follows that

\begin{equation}
W(z)J^a(w)=\frac{k+N}{(z-w)^2}W^a(z),
\end{equation}

take a derivative in  $w$, one obtains

\begin{equation}\label{h1}
W(z)\partial^r
J^a(w)=(-1)^r\frac{(k+N)\frac{(r+2)!}{2}}{(z-w)^{2+r}}W^a(z).
\end{equation}

From \eqref{h1} one gets the operator expansion

\begin{align}
\begin{split}
\label{oper11}
&W(z)\partial^r J^a(w)=(-1)^r(k+N)\frac{(r+2)!}{2}(\frac{W^a(w)}{(z-w)^{2+r}}+...+\frac{\frac{1}{k!}\partial^{k}W^a(w)}{(z-w)^{2+r-k}})+...\\
&+\frac{\frac{1}{(2+r-1)!}\partial^{2+r-1}W^a(w)}{(z-w)}).
\end{split}
\end{align}

%в частности, для $n=1$ имеем

%\begin{equation}\label{oper1}
%W(z)\partial^n J^a(w)=(k+N)(\frac{W^a(z)}{(z-w)^{2}}+\frac{\partial W^a(z)}{(z-%w)})+...,
%\end{equation}

Compare the formulas \eqref{f0} and \eqref{oper11},  one obtains the
theorem

\begin{thm}

The following formulas take place
\begin{align}
\begin{split}
% & Q_{-3}J^a(z)=0\\
%& Q_{-1}J^a(z)=(k+N)Q^a(z)\\
%& Q_{-2}J^a(z)=(k+N)\partial Q^a(z)\\
& W_{n}(\partial^r J^a(w))=0,\,\,n> r-1,\\
& W_{n}(\partial^r
J^a(w))=(-1)^r(k+N)\frac{(r+2)!}{2}\frac{1}{(r-n-1)!}\partial^{r-n-1}
W^a(w),\,\, n\leq r-1.
\end{split}
\end{align}

\end{thm}

Thus under the action of  $W_n$  the fields $J^a(w)$  and their
derivatives are mapped to  $W^a(w)$ and their derivatives.

From the formula \eqref{r2} one obtains that

\begin{equation}\label{h2}
W(z)\partial^r
W^{\gamma}(w)=(-1)^r(k+N)(k+\frac{N}{2})\frac{2}{N}(N^2-4)\frac{J^{\gamma}(z)}{(z-w)^{4+r}}.
\end{equation}

From \eqref{h2} one obtains the operator expansion

\begin{align}
\begin{split}
\label{oper12} & W(z)\partial^r
W^{\gamma}(w)=(-1)^r(k+N)(k+\frac{N}{2})\frac{2}{N}(N^2-4)
(\frac{J^{\gamma}(w)}{(z-w)^{4+r}}+...\\
&+\frac{\frac{1}{k!}\partial^kJ^{\gamma}(w)}{(z-w)^{4+r-k}}+..+
\frac{\frac{1}{(4+r-1)!}\partial^{4+r-1}J^{\gamma}(w)}{(z-w)^{}}).
\end{split}
\end{align}

Compare the formulas \eqref{f0} and \eqref{oper12}, one obtains the
theorem.

\begin{thm}

The following formulas take place
\begin{align}
\begin{split}
&W_n\partial^r W^a(w)=0, \,\,\, n>r+1\\
& W_n\partial^r
W^a(w)=(-1)^r(k+N)(k+\frac{N}{2})\frac{2}{N}(N^2-4)\frac{1}{(r-n+1)!}\partial^{r-n+1}J^{\gamma}(w).
%&Q_{1}Q^a(w)=(3k-N)(2k-N)\frac{2}{N}(N^2-4)2J^a(w)\\
%&Q_{0}Q^a(w)=(3k-N)(2k-N)\frac{2}{N}(N^2-4)\partial J^a(w)\\
%&Q_{-1}Q^a(w)=(3k-N)(2k-N)\frac{2}{N}(N^2-4)\partial^2 J^a(w)\\
%&Q_{-2}Q^a(w)=(3k-N)(2k-N)\frac{2}{N}(N^2-4)\partial^3 J^a(w)
\end{split}
\end{align}

\end{thm}

Thus under the action of  $W_n$  the fields $W^a(w)$  and their
derivatives are mapped to  $J^a(w)$ and their derivatives.

\begin{cor}
The operators  $W_n$ are not differential operators on the
$WZW$-primary fields, those are only the compositions $W_nW_m$.
\end{cor}

Thus it is shown explicitly that there exist no higher analogue of
the Knizhnik-Zamolodchikov equation associated with the algebra
$WA_2$ (that is with the Casimir element of the third order).
% The reason is that one can represent $W_n$ as an algebraic operator, but one can represent as
%  a differential operator only the composition $W_nW_m$ of two generators of the $W$-algebra
%  \footnote{Thus one can construct some equation correponding to compositions of generators of $W$-algebra}.

\section{The energy-momentum tensor associated with the Capelli element of the fourth order} \label{kzs2}

Let us construct an analog of the Knizhnik-Zamolodchikoc equation
using the Capelli element $C_4$.
\subsection{The correlation functions under consideration}
\label{corf}

Let us change  the type of considered correlation functions. for
this purpose let us define a field associated with a pfaffian

\begin{equation}
Pf F_{I}(z)=\frac{1}{k!2^k}\sum_{\sigma\in
S_{2k}}(-1)^{\sigma}(F_{\sigma(i_1)\sigma(i_2)}(F_{\sigma(i_3)\sigma(i_4)}(F_{\sigma(i_5)\sigma(i_6)}...F_{\sigma(i_{2k-1})\sigma(i_{2k})}))...)(z),
\end{equation}

 The classical Knizhnik-Zamolodchikov equation is system of PDE for
 the correlation function if type

\begin{equation}\label{cf1}
<\varphi_1(z_1)...\varphi_r(z_r)>,
\end{equation}

where $\varphi_k(z_k)$ are $WZW$ primary fields.

The equations that are constructed below are equations for the
correlation functions of the following type.
%As it was pointed in
%Section  \ref{prim}  a current $F_{ij}(z)$ in the case $k=0$ is a
% $WZW$-primary field that transform in the adjoint equation.
For the set of indices $I=\{i_1,i_2,i_3,i_{4}\}\subset \{1,...,N\}$
consider a field $\varphi(z)=PfF_I(z)$. This fields is not
$WZW$-primary even for $k=0$. Introduce a correlation function

\begin{equation}\label{cf2}
<PfF_I(z_1)\varphi_2(z_2)...\varphi_r(z_r)>,
\end{equation}

where $\varphi_k(z_k)$ are $WZW$-primary fields. Below an equation
for this correlation function is derived.

This correlation function is a function of variables $z_1,...,z_r$
that takes values in a $r$-th tensor power of an adjoint
representation of the algebra $\mathfrak{o}_N$.

\subsection{A field associated with a Capelli element of the order
$4$. The algebra $WB_2$}\label{wb}

Define fields

\begin{defn}

\begin{equation}
C_n(z)=\sum_{|I|=n}(PfF_IPfF_I)(z).
\end{equation}

\end{defn}

The algebra generated by $1,L_n,C_4^k$, is a Casimir $W$-algebra, it
is denoted $WB_2$.

For an arbitrary field $\varphi(w)$ one has

\begin{equation}
C_4(z)\varphi(w)=\sum_{k}(C_4^k\varphi)(w)(z-w)^{-k-4}.
\end{equation}

Consider an action of  $C_4^{-1}$ on the field $PfF_I(z_1)$. From
one hand let us represent this action as a differential operator.
More precise we show that

\begin{equation}
C_4^{-1}PfF_I(z_1)=const\partial PfF_I(z_1)
\end{equation}

From the other hand we show that the correlation functions

\begin{equation}
<C_4^{-1}PfF_I(z_1)\varphi_2(z_2)...\varphi_r(z_r)>
\end{equation}

can be algebraically expressed through the correlation functions of
type \eqref{cf1}, \eqref{cf2}. When one equates two these
expressions one obtains a higher analog of the
Knizhnik-Zamolodchikov equation.

\subsection{Representation of  $C_4^{-1}$ as a differential operator}

To represent $C_4^{-1}$  as a differential operator it is necessary
to find two  higher  terms of the OPE $C_4(z)PfF_J(w)$.

As a first step let us calculate the contraction
$\overbrace{PfF_J(z)PfF_I(x)}$.

 \subsubsection{Calculation of $J^a(z)(J^{\beta} J^{\gamma})(w)$}

Let us do the following general calculation. Let $J^{\alpha}$ be an
orthonomal with respect to the Killing form base of a simple Lie
algebra. Let $f_{\alpha,\beta,\gamma}$ be structure constants.

Calculate the contraction

$$\overbrace{J^a(z)(J^{\beta} }J^{\gamma})(w).$$

The Proposition takes place

\begin{prop}
\label{p1}
\begin{align}
\begin{split}
& \overbrace{J^a(z)(J^{\beta}
}J^{\gamma})(w)=\frac{k\delta_{a,\beta} J^{\gamma}(w)}{(z-w)^2}+
\frac{k\delta_{a,\gamma}J^{\beta}(w)}{(z-w)^2}+
\frac{f_{a,\beta,\gamma}k}{(z-w)^3}+\\&+\frac{f_{a,\beta,c}f_{c,\gamma,d}J^d(w)}{(z-w)^2}+
\frac{f_{a,\beta,c}(J^cJ^{\gamma})(w)}{z-w}+\frac{f_{a,\gamma,c}(J^{\beta}J^c)(w)}{z-w}
\end{split}
\end{align}

\end{prop}

The proof of this Proposition can be found in Appendix in Section
\ref{ap1}

\subsubsection{Calculation of
$F_{j_1,j_2}(z)(F_{i_1,i_2}F_{i_3,i_4})(w)$}

In the case of the orthogonal algebra the formula from the
Proposition \ref{p1} looks as follows.

\begin{prop}
\label{p2}
\begin{align}
\begin{split}
&
\overbrace{F_{j_1,j_2}(z)(F_{i_1,i_2}}F_{i_3,i_4})(w)=\frac{k\delta_{(j_1,j_2),(i_1,i_2)}F_{i_3,i_4}}{(z-w)^2}+\frac{k\delta_{(j_1,j_2),(i_3,i_4)}F_{i_1,i_2}}{(z-w)^2}+\\
&+\frac{F_{F_{i_1,i_2}F_{i_3,i_4}\{j_1,j_2\}}(w)}{(z-w)^2}+\frac{F_{F_{j_1,j_2}\{i_1,i_2\}}F_{i_3,i_4}(w)}{(z-w)^2}+\frac{F_{i_1,i_2}F_{F_{j_1,j_2}\{i_3,i_4\}}}{(z-w)}.
\end{split}
\end{align}
\end{prop}

\subsubsection{Calculation of $F_{j_1,j_2}(z)PfF_I(w)$  in the case $|I|$=4}

Let us find the contraction $\overbrace{F_{j_1,j_2}(z)(PfF_I)(w)}$
in the case $|I|=4$.

 Let $$I=\{i_1,i_2,i_3,i_4\}.$$

One has

\begin{align}
\begin{split}
& \overbrace{F_{j_1,j_2}(z)(PfF_I)(w)}=
\frac{k(-1)^{\{j_1,j_2\},I\setminus\{j_1,j_2\}}F_{I\setminus\{j_1,j_2\}}(w)}{(z-w)^2}+\\&
 +\frac{F_{PfF_I\{j_1,j_2\}}(w)}{(z-w)^2}+\frac{(PfF_{F_{j_1,j_2}I})(w)}{(z-w)}.
 \end{split}
\end{align}

After simplification one gets

\begin{prop}
\label{p4}
\begin{align}
\begin{split}
& \overbrace{F_{j_1,j_2}(z)(PfF_I)(w)}=
\frac{(k+2)F_{PfF_I\{j_1,j_2\}}(w)}{(z-w)^2}+\frac{(PfF_{F_{j_1,j_2}I})(w)}{(z-w)}.
 \end{split}
\end{align}
\end{prop}

\subsubsection{ The calculation of a contraction of two pfaffians
   }

Let $|I|=|J|=4$. One has

\begin{prop}
\label{p5}
\begin{align}
\begin{split}
 &\overbrace{PfF_J(z)PfF_I(x)}=\frac{1}{2}\sum_{I=I_1\sqcup
 I_2,|I_1|=|I_2|=2}(-1)^{(I_1,I_2)}(
\frac{2(k+2)kPfF_{J\setminus I_1I_2}(z)}{(z-x)^4}-\\&-
\frac{2(k+2)PfF_{F_{I_2}(J\setminus I_1)}(z)}{(z-x)^3}+
\frac{(k+2)(F_{I_2}PfF_{J\setminus I_1})(z)}{(z-x)^2}-\\&-
\frac{(k+2)(F_{I_2}PfF_{J\setminus I_1})_1(z)}{(z-x)}-
\frac{(F_{I_2}PfF_{F_{I_1}J})(z)}{(z-x)}+\\&+\frac{(k+2)(F_{I_1}PfF_{J\setminus
I_2})(z)}{(z-x)^2}+\frac{(k+2)(F_{I_1}PfF_{J\setminus
I_2})_1(z)}{(z-x)}- \frac{(F_{I_1}PfF_{F_{I_2}J})(z)}{(z-x)}),
\end{split}
\end{align}

where $(-1)^{(I_1,I_2)}$ is a sign of a permutation of the set $I$,
that firstly puts the set $I_1$, and then the set $I_2$.
\end{prop}

\proof

The following formula  takes place  which is a particular case of
the minor summation formula \cite{AG}

\begin{align}
\begin{split}
PfF_I=\frac{1}{2}\sum_{I=I_1\sqcup
 I_2,|I_1|=|I_2|=2}(-1)^{(I_1,I_2)} F_{I_1}F_{I_2}.
\end{split}
\end{align}

Using it one gets

\begin{align}
\begin{split}
 &\overbrace{PfF_J(z)PfF_I(x)}=\frac{1}{2}\sum_{I=I_1\sqcup
 I_2,|I_1|=|I_2|=2}(-1)^{(I_1,I_2)}
 \frac{1}{2\pi i}\oint \frac{dx_1}{(x_1-x)}
\overbrace{PfF_J(z)F_{I_1}(x_1)}F_{I_2}(x)+\\&+F_{I_1}(x_1)\overbrace{PfF_J(z)F_{I_2}(x)}.
 \end{split}
\end{align}

Substitute this into this expression the formula for the contraction
of a pfaffian and a current that was found in Proposition \ref{p4}.
One gets

\begin{align}
\begin{split}
&\frac{1}{2}\sum_{I=I_1\sqcup
 I_2,|I_1|=|I_2|=2}(-1)^{(I_1,I_2)}\frac{1}{2\pi
i}\oint \frac{dx_1}{(x_1-x)} (\frac{(k+2)PfF_{J\setminus
I_1}(z)}{(z-x_1)^2}-\\&-\frac{PfF_{F_{I_1}J}(z)}{(z-x_1)})F_{I_2}(x)+F_{I_1}(x_1)(\frac{(k+2)PfF_{J\setminus
I_2}(z)}{(z-x)^2}-\frac{PfF_{F_{I_2}J}(z)}{(z-x)})=\\ &=
\frac{1}{2}\sum_{I=I_1\sqcup
 I_2,|I_1|=|I_2|=2}(-1)^{(I_1,I_2)}\frac{1}{2\pi i}\oint
\frac{dx_1}{(x-x_1)}(\frac{(k+2)}{(z-x_1)^2}PfF_{J\setminus
I_1}(z)F_{I_2}(x)-\\&-\frac{1}{(z-x_1)}PfF_{F_{I_1}J}(z)F_{I_2}(x)+
\frac{(k+2)}{(z-x)^2}F_{I_1}(x_1)PfF_{J\setminus
I_2}(z)-\frac{1}{(z-x)}F_{I_1}(x_1)PfF_{F_{I_2}J}(z)).
 \end{split}
\end{align}

After integration and application of the formula from Proposition
\ref{p4} one gets

\begin{align}
\begin{split}
\label{long1}& \frac{1}{2}\sum_{I=I_1\sqcup
 I_2,|I_1|=|I_2|=2}(-1)^{(I_1,I_2)}(
\frac{(k+2)kPfF_{J\setminus I_1I_2}(z)}{(z-x)^4}-
\frac{(k+2)PfF_{F_{I_2}(J\setminus I_1)}(z)}{(z-x)^3}+\\&+
\frac{(k+2)(F_{I_2}PfF_{J\setminus I_1})(z)}{(z-x)^2}-
\frac{(k+2)(F_{I_2}PfF_{J\setminus I_1})_1(z)}{(z-x)}-
\frac{(k+2)PfF_{(F_{I_1}J)\setminus I_2}(z)}{(z-x)^3}+\\&+
\frac{PfF_{F_{I_1}F_{I_2}J}(z)}{(z-x)^2}-
\frac{(F_{I_2}PfF_{F_{I_1}J})(z)}{(z-x)}+\frac{(k+2)kPfF_{J\setminus
I_1I_2}(z)}{(z-x)^4}+\\&+\frac{(k+2)PfF_{F_{I_1}(J\setminus
I_2)}(z)}{(z-x)^3}+\frac{(k+2)(F_{I_1}PfF_{J\setminus
I_2})(z)}{(z-x)^2}+\\&+\frac{(k+2)(F_{I_1}PfF_{J\setminus
I_2})_1(z)}{(z-x)}-\frac{(k+2)PfF_{F_{I_2}J\setminus
I_1}(z)}{(z-x)^3}- \frac{PfF_{F_{I_1}F_{I_2}}(z)}{(z-x)^2}-
\frac{(F_{I_1}PfF_{F_{I_2}J})(z)}{(z-x)})
 \end{split}
\end{align}

After simplification one gets

\begin{align}
\begin{split}
\label{long}& \frac{1}{2}\sum_{I=I_1\sqcup
 I_2,|I_1|=|I_2|=2}(-1)^{(I_1,I_2)}(
\frac{2(k+2)kPfF_{J\setminus I_1I_2}(z)}{(z-x)^4}-
\frac{2(k+2)PfF_{F_{I_2}(J\setminus I_1)}(z)}{(z-x)^3}+\\&+
\frac{(k+2)(F_{I_2}PfF_{J\setminus I_1})(z)}{(z-x)^2}-
\frac{(k+2)(F_{I_2}PfF_{J\setminus I_1})_1(z)}{(z-x)}-
\frac{(F_{I_2}PfF_{F_{I_1}J})(z)}{(z-x)}+\\&+\frac{(k+2)(F_{I_1}PfF_{J\setminus
I_2})(z)}{(z-x)^2}+\frac{(k+2)(F_{I_1}PfF_{J\setminus
I_2})_1(z)}{(z-x)}- \frac{(F_{I_1}PfF_{F_{I_2}J})(z)}{(z-x)})
 \end{split}
\end{align}
\endproof

\subsubsection{Calculation of the contraction of $PfF_J(z)$ and $C_4(w)$}

Let us prove the Theorem

\begin{thm}
\begin{align}
\begin{split}
 &\overbrace{PfF_J(z)С_4(w)}=(6(k+2)k+12(N-4)-(k+2)((N-2)(N-3)(N-4)+6))\frac{PfF_J(w)}{(z-w)^4}+\\&+
 (12(N-4)-(k+2)((N-2)(N-3)(N-4)+6))\frac{4\partial PfF_J(w)}{(z-w)^3}+l.o.t.\footnote{l.o.t. = lower order terms}
 \end{split}
\end{align}

\end{thm}
\proof One has

\begin{align}
\begin{split}
\label{kuda}
 &\overbrace{PfF_J(z)(PfF_I}PfF_I)(w)=\\&=\frac{1}{2\pi
i}\oint
\frac{dx}{x-w}\overbrace{PfF_J(z)PfF_I(x)}PfF_I(w)+PfF_I(x)\overbrace{PfF_J(z)PfF_I(w)}.
 \end{split}
\end{align}

This expression consists of two terms. Let us show that they are
equal. Thus the expression \eqref{kuda} equals to

\begin{equation}\label{t1}
2\frac{1}{2\pi i}\oint
\frac{dx}{x-w}\overbrace{PfF_J(z)PfF_I(x)}PfF_I(w).
\end{equation}

Indeed, let us introduce notations $a(z)=PfF_J(z)$, $b(x)=PfF_I(x)$,
$c(w)=PfF_I(w)$. Proposition \ref{p5} gives us an expression of type

\begin{equation}
a(z)b(w)=\sum_m (ab)_{m}(z)(z-w)^m.
\end{equation}

Then \eqref{kuda} can be written as follows

\begin{align}
\begin{split}
\label{kuda1}
 &\frac{1}{2\pi
i}\oint
\frac{dx}{x-w}\overbrace{a(z)b(x)}c(w)+c(x)\overbrace{a(z)b(w)}=\\
&=\sum_m (ab)_m(z)c(w)(z-w)^m+c(w)(ab)_m(z)(z-w)^m.
 \end{split}
\end{align}

However  $(ab)_m(z)c(w)(z-w)^m=c(w)(ab)_m(z-w)^m$, that is why the
two summands in \eqref{kuda} are equal.

Consider now the expression \eqref{t1}. Let us prove that when one t
calculates contraction  in the numerator there appear  no singular
terms  in \eqref{t1} of order greater then $(z-w)^{-4}$, the term of
order $(z-w)^{-4}$ is proportional $PfF_J(w)(z-w)^{-4}$, the term of
order $(z-w)^{-3}$ is proportional to $ PfF_J(w)(z-w)^{-3}$.

To prove this let us substitute the expression \eqref{long} into
\eqref{t1} and then consider the terms separately. Thus it is
necessary to investigate the following expressions.
% Since the
%summation over  $I_1,I_2$ is taken it is necessary to consider only
%the following terms form \eqref{long}

\begin{enumerate}
\item $\sum_{I}\sum_{I=I_1\sqcup I_2}(-1)^{(I_1,I_2)}\frac{1}{2\pi i}\oint \frac{dx}{x-w}\frac{PfF_{J\setminus I_1I_2}(z)}{(z-x)^4}PfF_I(w)$,
\item $\sum_{I}\sum_{I=I_1\sqcup I_2}(-1)^{(I_1,I_2)}\frac{1}{2\pi i}\oint \frac{dx}{x-w}\frac{PfF_{F_{I_2}(J\setminus I_1)}(z)}{(z-x)^3}PfF_I(w)$,
\item $\sum_{I}\sum_{I=I_1\sqcup I_2}(-1)^{(I_1,I_2)}\frac{1}{2\pi i}\oint \frac{dx}{x-w}\frac{(F_{I_2}PfF_{J\setminus I_1})(z)}{(z-x)^2}PfF_I(w)$,
\item $\sum_{I}\sum_{I=I_1\sqcup I_2}(-1)^{(I_1,I_2)}\frac{1}{2\pi i}\oint \frac{dx}{x-w}\frac{(F_{I_2}PfF_{J\setminus I_1})_1(z)}{(z-x)}PfF_I(w)$,
\item $\sum_{I}\sum_{I=I_1\sqcup I_2}(-1)^{(I_1,I_2)}\frac{1}{2\pi i}\oint \frac{dx}{x-w}\frac{(F_{I_1}PfF_{F_{I_2}J})(z)}{(z-x)}PfF_I(w)$.
%Вычисляем почленно с членами $\overbrace{PfF_J(z)PfF_I(x)}$.
\end{enumerate}

The summation $\sum_{I}$ is taken over all subsets $I\subset \{1,...,N\}$ that
 consist of four elements and the summation $\sum_{I=I_1\sqcup I_2}$ is taken over all partitions of
 $I$  into two subsets that consist of two elements\footnote{Everywhere below we shall use these notations.}.

Consider  cases  in further five lemmas.

\begin{lem}
\label{l1}
\begin{align}
\begin{split}
& \sum_{I}\sum_{I=I_1\sqcup I_2} \frac{1}{2\pi i}\oint
\frac{dx}{x-w}\frac{(k+2)kPfF_{J\setminus
I_1I_2}(z)}{(z-x)^4}Pf_I(w)=\frac{6(k+2)kPfF_J(w)}{(z-w)^4}.
 \end{split}
\end{align}
\end{lem}

This Lemma is obvious.

%\begin{align}
%\begin{split}
%& \frac{1}{2}\sum_{I=I_1\sqcup
% I_2,|I_1|=|I_2|=2}(-1)^{(I_1,I_2)} \frac{1}{2\pi i}\oint \frac{dx}{x-w}\frac{PfF_{J\setminus
%I_1I_2}(z)}{(z-x)^4}Pf_I(w)=\frac{PfF_J(w)}{2(z-w)^4}.
% \end{split}
%\end{align}

\begin{lem}
\label{l2}
\begin{align}
\begin{split}
&  \sum_{I}\sum_{I=I_1\sqcup I_2} (-1)^{(I_1,I_2)}\frac{1}{2\pi
i}\oint \frac{dx}{x-w} \frac{PfF_{F_{I_2}(J\setminus
I_1)}(z)}{(z-x)^3}PfF_I(w)=\\&=(N-4)(\frac{12PfF_J(w)}{(z-w)^{4}}+\frac{2\partial
PfF_J(w)}{(z-w)^3}+l.o.t.)
\end{split}
\end{align}
\end{lem}

The proof of the Lemma \ref{l2} can be found in Application
\ref{ap4}.

\begin{lem}
\label{l3} In the OPE \begin{align}
\begin{split}\label{64}
&\sum_{I}\sum_{I=I_1\sqcup I_2}(-1)^{(I_1,I_2)}\frac{1}{2\pi i}\oint
\frac{dx}{x-w}\frac{(F_{I_2}PfF_{J\setminus
I_1})(z)}{(z-x)^2}PfF_I(w)
\end{split}
\end{align}
there are no singular terms of order greater than $(z-w)^{-2}$.
\end{lem}

The proof of this Lemma can be found in Application  \ref{ap5}

%{\bf 5. The case } $\frac{(F_{I_1}PfF_{F_{I_2}J})(z)}{(z-w)}$.

\begin{lem}
\label{l4} In the OPE \begin{align}
\begin{split}\label{644}
&\sum_{I}\sum_{I=I_1\sqcup I_2}(-1)^{(I_1,I_2)}\frac{1}{2\pi i}\oint
\frac{dx}{x-w}\frac{(F_{I_2}PfF_{J\setminus
I_1})_1(z)}{(z-x)^2}PfF_I(w)
\end{split}
\end{align}
there are no singular terms of order greater than $(z-w)^{-2}$.
\end{lem}

This lemma is proved analogously to Lemma \ref{64}.

\begin{lem}
\label{l5}
\begin{align}
\begin{split}\label{167}
&\sum_{I}\sum_{I=I_1\sqcup I_2}(-1)^{(I_1,I_2)} \frac{1}{2\pi
i}\oint \frac{dx}{x-w}
 \frac{(F_{I_1}PfF_{F_{I_2}J})(z)PfF_I(w)}{(z-x)}=\\
 &=-(k+2)((N-2)(N-3)(N-4)+6)(\frac{PfF_I(w)}{(z-w)^4}+\frac{\partial
 PfF_I(w)}{(z-w)^3}+l.o.t.)
\end{split}
\end{align}
\end{lem}

This Lemma is proved in Application  \ref{ap6}.

We have obtained that

\begin{align}
\begin{split}
 &\overbrace{PfF_J(z)С_4(w)}=(6(k+2)k+12(N-4)-(k+2)((N-2)(N-3)(N-4)+6))\frac{PfF_J(w)}{(z-w)^4}+\\&+
 (12(N-4)-(k+2)((N-2)(N-3)(N-4)+6))\frac{\partial PfF_J(w)}{(z-w)^3}+l.o.t.
 \end{split}
\end{align}

The theorem is proved.

\endproof

If one  changes the fields one obtains the OPE

\begin{align}
\begin{split}
 &\overbrace{С_4(z)PfF_J(w)}=(6(k+2)k+12(N-4)+(k+2)((N-2)(N-3)(N-4)-6))\frac{PfF_J(w)}{(z-w)^4}+\\&+
 (6(k+2)k)\frac{\partial PfF_J(w)}{(z-w)^3}+l.o.t.
 \end{split}
\end{align}

Thus one has

\begin{cor}\label{s1}
\begin{equation}
C_4^{-1}PfF_J(w)=6(k+2)k\partial PfF_J(w).
\end{equation}
\end{cor}

\subsection{Representation of  $C_4^{-1}$ as an algebraic operator}

In this Section an algebraic operator   $\mathcal{A}$ is constructed
that expresses the correlation function

\begin{equation}
<C_4^{-1}PfF_I(z_1)\varphi_2(z_2)...\varphi_r(z_r)>,
\end{equation}

 through the correlation functions of type \eqref{cf1},
\eqref{cf2}.

Consider the field  $PfF_I(z)$.  For it the following OPE takes
place

\begin{align}
\begin{split}
& \overbrace{F_{J}(z)PfF_I(w)}=\sum_m
(F_{J}^kPfF_I)(w)(z-w)^{-m-1}=\\&=\frac{(k+2)PfF_{I\setminus
J}(w)}{(z-w)^2}+\frac{PFF_{F_J}I(w)}{(z-w)}.
 \end{split}
\end{align}

From this equality one gets

\begin{align}
\begin{split}\label{jpf}
& F_J^{n}PfF_I(w)=0,\,\,\, n>1,\\
& F_J^{1}PfF_I(w)=(k+2)PfF_{I\setminus J}(w),\\
& F_J^0PfF_I(w)=PFF_{F_J}I(w).
 \end{split}
\end{align}

Let use the general fact

 \begin{equation}
(AB)_m=\sum_{n\leq h_A}A_nB_{m-n}+\sum_{n>h_A}B_{m-n}A_n.
 \end{equation}

Then one has

\begin{align}
\begin{split}
&C_4^{-1}=(PfF_J)_{2}(PfF_J)_{-3}+(PfF_J)_{1}(PfF_J)_{-2}+...\\&+(PfF_J)_{-4}(PfF_J)_{3}+(PfF_J)_{-5}(PfF_J)_{4}+...,
 \end{split}
\end{align}

where

\begin{equation}
(PfF_J)_m=\frac{1}{2}\sum_{J=J_1\sqcup
J_2}(-1)^{(J_1,J_2)}F_{J_1}^{1}F_{J_2}^{m-1}+F_{J_1}^{0}F_{J_2}^{m}+...+F_{J_2}^{m-2}F_{J_2}^{2}+F_{J_2}^{m-3}F_{J_1}^{3}+...
\end{equation}

 From here one gets

\begin{align}
\begin{split}\label{rq2}
& (PfF_J)_mPfF_I=0,  \,\,\, m>2.
\end{split}
\end{align}

Thus using \eqref{jpf}, \eqref{rq2}, one gets that

\begin{align}
\begin{split}\label{alc}
& C_4^{-1}PfF_I(w)=\frac{1}{4}\sum_{J}\sum_{J=J_1\sqcup
J_2}\sum_{J=J'_1\sqcup J'_2}\sum_{k,l,p,q}
(-1)^{(J_1,J_2)}(-1)^{(J'_1,J'_2)}F_{J_1}^kF_{J_2}^l F_{J'_1}^p
F_{J'_2}^qPfF_I(w),
 \end{split}
\end{align}
where the following sets of indices $(k,l,p,q)$ are admissible

%\begin{align}
%\begin{split}
%&(1,1,0,-3),(1,1,1,-4),(1,1,-1,-2)\\
%& (0,1,-1,-1),(0,1,0,-2),(0,1,1,-3)\\
%& (0,0,0,-1),(0,0,1,-2),(1,-1,0,-1),(1,-1,1,-2)\\
%& (0,-1,0,0),(0,-1,1,-1),(1,-2,0,0),(1,-2,0,0),(1,-2,1,-1)\\
%& (-1,-1,0,1),(1,-3,0,1)\\
%& (0,-3,1,1), (1,-4,1,1)
% \end{split}
%\end{align}

\begin{align}
\begin{split}
& (0,1,-1,-1),(0,1,0,-2), (0,0,0,-1),
 \end{split}
\end{align}

and also all other collections that are obtained from these by any
permutation of elements.

Let us prove the theorem

\begin{thm}
\label{ta}
\begin{align}
\begin{split}
&<C_4^{-1}PfF_I(z_1)\varphi_2(z_2)...\varphi_r(z_r)>
=\\&=\mathcal{A}(<PfF_I(z_1)\varphi_2(z_2)...\varphi_r(z_r)>,<\varphi(z_1)\varphi_2(z_2)...\varphi_r(z_r)>),
 \end{split}
 \end{align}
 where the operator $\mathcal{A}$ is algebraic.

\end{thm}

\proof Let us use the formula \eqref{alc}. Prove that the summand
corresponding to every collection of indices $(k,l,p,q)$
  is represented as  a result of an application of an algebraic operator  to
  correlation fucntions of type \eqref{cf1},
\eqref{cf2}. Let us write explicitely the formulas for the action of
operators  $F_J^{-2}$, $F_J^{-1}$, $F_J^0$, $F_J^1$.

Let us consider first the simplest case. One has
\begin{align}
\begin{split}
<F_{J}^0PfF_I(z_1)\varphi_2(z_2)...\varphi_r(z_r)>=F_J^{(1)}<PfF_I(z_1)\varphi_2(z_2)...\varphi_r(z_r)>,
\end{split}
\end{align}

where $F_J^{(k)}$ is the operator $F_J$,  acting on the $k$-th
tensor power (the correlation function takes values in the $m$-th
tensor power of the adjoint representation of $\mathfrak{o}_N$ - see
Section \ref{corf}).

Standard calculation in the $WZW$ theory show that

\begin{align}
\begin{split}
<F_{J}^{-1}PfF_I(z_1)\varphi_2(z_2)...\varphi_r(z_r)>=-\sum_{j\geq
2}\frac{F_J^{(j)}}{(z_1-z_j)}<PfF_I(z_1)\varphi_2(z_2)...\varphi_r(z_r)>.
\end{split}
\end{align}

Analogously one gets

\begin{align}
\begin{split}
<F_{J}^{-2}PfF_I(z_1)\varphi_2(z_2)...\varphi_r(z_r)>=-\sum_{j\geq
2}\frac{F_J^{(j)}}{(z_1-z_j)^2}<PfF_I(z_1)\varphi_2(z_2)...\varphi_r(z_r)>.
\end{split}
\end{align}

Finally  by \eqref{jpf} one gets

\begin{align}
\begin{split}
<F_{J}^{1}PfF_I(z_1)\varphi_2(z_2)...\varphi_r(z_r)>=(k+2)<PfF_{I\setminus
J}(z_1)\varphi_2(z_2)...\varphi_r(z_r)>,
\end{split}
\end{align}

note that the field $PfF_{I\setminus J}(z_1)=F_{I\setminus J}(z_1)$
is not primary.

\section{Higher Knizhnik-Zamolodchikov equations}
\label{ur} Using the Corollary \ref{s1}  and Theorem \ref{ta} one
obtains the equation

\begin{align}
\begin{split}
&6(k+2)k\partial_{z_1}<PfF_J(z_1)\varphi_2(z_2)...\varphi_r(z_r)>=\\
&=\frac{1}{4}\sum_{I}\sum_{I=I_1\sqcup I_2,I'=I'_1\sqcup I'_2}
(-1)^{(I_1,I_2)}(-1)^{(I'_1,I'_2)}\\
&2\sum_{j\geq 2}\frac{F_{I_1}^{(1)}F_{I_2}^{(1)} F_{I'_1}^{(1)}
F_{I'_2}^{(j)}+F_{I_1}^{(1)}F_{I_2}^{(j)} F_{I'_1}^{(1)}
F_{I'_2}^{(1)}}{(z_1-z_j)}<PfF_J(z_1)\varphi_2(z_2)...\varphi_r(z_r)>+\\
&+(k+2)(8\sum_{i\leq j\geq 2}\frac{F_{I_2}^{(1)} F_{I'_1}^{(i)}
F_{I'_2}^{(j)}+F_{I_2}^{(i)} F_{I'_1}^{(1)}
 F_{I'_2}^{(j)}+F_{I_2}^{(i)} F_{I'_1}^{(j)} F_{I'_2}^{(1)}}{(z_i-z_1)(z_j-z_1)}-\\
&-4\sum_{j\geq 2}\frac{F_{I_2}^{(1)} F_{I'_1}^{(1)}
F_{I'_2}^{(j)}+F_{I_2}^{(1)} F_{I'_1}^{(1)}
F_{I'_2}^{(j)}+F_{I_2}^{(j)} F_{I'_1}^{(1)}
F_{I'_2}^{(1)}}{(z_j-z_1)^2})<PfF_{J\setminus
I_1}(z_1)\varphi_2(z_2)...\varphi_r(z_r)>+\\
&+(k+2)(8\sum_{i\leq j\geq 2}\frac{F_{I_1}^{(1)} F_{I_2}^{(i)}
F_{I'_2}^{(j)}+F_{I_1}^{(i)}
 F_{I_2}^{(1)} F_{I'_2}^{(j)}+F_{I_1}^{(i)} F_{I_2}^{(j)} F_{I'_2}^{(1)}}{(z_i-z_1)(z_j-z_1)}-\\
&-4\sum_{j\geq 2}\frac{F_{I_1}^{(1)} F_{I_2}^{(1)}
F_{I'_2}^{(j)}+F_{I_1}^{(1)} F_{I_2}^{(1)}
F_{I'_2}^{(j)}+F_{I_1}^{(j)} F_{I_2}^{(1)}
F_{I'_2}^{(1)}}{(z_j-z_1)^2})<PfF_{J\setminus
I'_1}(z_1)\varphi_2(z_2)...\varphi_r(z_r)>
\end{split}
\end{align}

The summation $\sum_I$ is taken over all subsets $I\subset
\{1,...,N\}$ that consist of four elements. The summation
$\sum_{I=I_1\sqcup I_2,I'=I'_1\sqcup I'_2}$ is taken over all
partitions $I=I_1\sqcup I_2$, $I'=I'_1\sqcup I'_2$ of the set  $I$
into two two-element subsets. The operator  $F_{I}^{(k)}$ is the
operator $F_I$, acting onto the $k$-th tensor component of the
correlation function.

To simplify the equation we use the fact that for $I_1\cap
I_2=\emptyset$, $I'_1\cap I'_2=\emptyset$, one has
$[F^{(k)}_{I_1},F_{I_2}^{(l)}]=[F^{(k)}_{I'_1},F_{I'_2}^{(l)}]=0$.

Thus the sum that corresponds to the set $(0,0,0,-1)$ and it's
permutations can be written as follows

\begin{align}
\begin{split}
&\sum_{j\geq 2}F_{I_1}^{(1)}F_{I_2}^{(1)} F_{I'_1}^{(1)}
F_{I'_2}^{(j)}+F_{I_1}^{(1)}F_{I_2}^{(j)} F_{I'_1}^{(1)}
F_{I'_2}^{(1)}+F_{I_1}^{(1)}F_{I_2}^{(1)} F_{I'_1}^{(j)}
F_{I'_2}^{(1)}+F_{I_1}^{(j)}F_{I_2}^{(1)} F_{I'_1}^{(1)}
F_{I'_2}^{(1)}=\\&=2\sum_{j\geq 2}F_{I_1}^{(1)}F_{I_2}^{(1)}
F_{I'_1}^{(1)} F_{I'_2}^{(j)}+F_{I_1}^{(1)}F_{I_2}^{(j)}
F_{I'_1}^{(1)} F_{I'_2}^{(1)}.
\end{split}
\end{align}

Analogously one can simplify the sums to the sets $(0,1,0,-2))$,
$(0,1,-1,-1)$. As the result one obtains the equation written above.

\section{Conclusion}
We investigate the possibility of construction  of
Knizhnik-Zamolodchikov type in the $WZW$-theory using higher order
Casimir elements.

We consider the $W$-algebra that is generated by modes of the field
constructed from a Gelfand element of the third order. In this case
the construction of the equation is impossible.

Also we consider the $W$-algebra that is generated by modes of the
field constructed from a Capelli element. Here we manage to express
the action of the operator $W_{-1}$ though a differential operator.
Hence it is possible to construct a Knizhnik-Zamolodchikov type
equation.

Considerations  are done in the case of fields that transform in the
adjoint representation. Analogous construction in the case of
arbitrary tensor representation can be obtained by taking a tensor
power of the adjoint representation and passing to irreducible
components.

Note that there is a well-known relation between the
Knizhnik-Zamolodchikov equation and Gaudin hamiltonians. Gaudin
hamiltonians are right sides of Knizhnik-Zamolodchikov equations.
 At the same time there exist higher Gaudin hamiltonians
  that are actually higher Sugawara elements  \cite{FFR}. Higher Gaudin hamiltonians commute with
   Gaudin hamiltonians and thus have the same eigenfunctions. This indicates that
   the higher Gaudin hamiltonians are candidates for right sides of higher Knizhnik-Zamolodchikov equations.
   But actually it is impossible to obtain such equation. The reason is that it is impossible to express the
   action of the higher Gaudin hamiltonians as differential
   operators (see for example \cite{wrev}, \cite{r}).

\section{Appendix}\label{ap}

\subsection{Proof of Proposition \ref{p1}}

\begin{align}
\begin{split}
& \overbrace{J^a(z)(J^{\beta} }J^{\gamma})(w)=\frac{1}{2\pi i}\oint
\frac{dx}{x-w}\overbrace{J^a(z)J^{\beta}(x) }J^{\gamma}(w)+
J^{\beta}(x)\overbrace{J^a(z) J^{\gamma}(w)}.
\end{split}
\end{align}

The integration is taken along a contour that pass around $w$ the
point $z$ is not contained inside the contour.

One has

\begin{align}
\begin{split}
& J^{a}(z)J^{\beta}(x)=\frac{k\delta_{a,\beta}}{(z-x)^2}+\sum_c\frac{f_{a,\beta,c}J^{c}(x)}{(z-x)}\\
& J^{a}(z)J^{\gamma}(w)=\frac{k\delta_{a,\gamma}}{(z-w)^2}+\sum_c\frac{f_{a,\gamma,c}J^{c}(w)}{(z-w)}\\
\end{split}
\end{align}

Thus one gets

\begin{align}
\begin{split}
\label{jk1} &\frac{1}{2\pi i}\oint
\frac{dx}{x-w}\overbrace{J^a(z)J^{\beta}(x) }J^{\gamma}(w)+
J^{\beta}(x)\overbrace{J^a(z) J^{\gamma}(w)}=\\&=\frac{1}{2\pi
i}\oint
\frac{dx}{x-w}(\frac{k\delta_{a,\beta}}{(z-x)^2}+\sum_c\frac{f_{a,\beta,c}J^{c}(x)}{(z-x)}+)J^{\gamma}(w)+
J^{\beta}(x)(\frac{k\delta_{a,\gamma}}{(z-w)^2}+\sum_c\frac{f_{a,\gamma,c}J^{c}(w)}{(z-w)})
\end{split}
\end{align}

The following expansions take place

\begin{align}\begin{split}
&
\frac{1}{(z-x)^2}=\frac{1}{(z-w)-(x-w))^2}=\frac{1}{(z-w)^2}\frac{1}{(1-\frac{x-w}{z-w})^2}=\\
&=
\frac{1}{(z-w)^2}(1+\frac{(x-w)^1}{(z-w)^1}+\frac{(x-w)^2}{(z-w)^2}+\frac{(x-w)^3}{(z-w)^3}+...)^2=\\
&=
\frac{1}{(z-w)^2}(1+2\frac{(x-w)^1}{(z-w)^1}+3\frac{(x-w)^2}{(z-w)^2}+4\frac{(x-w)^3}{(z-w)^3}+...).
\end{split}\end{align}

Using these expansions one gets

\begin{align}
\begin{split}
\label{jk} &\frac{1}{2\pi i}\oint
\frac{dx}{x-w}(\frac{k\delta_{a,\beta}}{(z-x)^2}J^{\gamma}(w)+\frac{k\delta_{a,\gamma}}{(z-w)^2}J^{\beta}(x))=\\&=
\frac{k\delta_{a,\beta}}{(z-w)^2}J^{\gamma}(w)+\frac{k\delta_{a,\gamma}}{(z-w)^2}J^{\beta}(w).
\end{split}
\end{align}

After substitution of \eqref{jk} into  \eqref{jk1} one gets

\begin{align}
\begin{split}
\label{jk3}
 \frac{k\delta_{a,\beta} J^{\gamma}(w)}{(z-w)^2}+
\frac{k\delta_{a,\gamma}J^{\beta}(w)}{(z-w)^2}+\sum_c\frac{1}{2\pi
i}\oint\frac{dx}{x-w}\frac{f_{a,\beta,c}J^{c}(x)J^{\gamma}(w)}{(z-x)}+\frac{f_{a,\gamma,c}J^{\beta}(x)J^{c}(w)}{(z-w)}
\end{split}
\end{align}

One has the expansions

\begin{align}
\begin{split}
\label{jk4}
&J^{с}(x)J^{\gamma}(w)=\frac{k\delta_{c,\gamma}}{(x-w)^2}+\sum_c\frac{f_{c,\gamma,d}J^{d}(w)}{(x-w)}+(J^{c}J^{\gamma})(w)+...\\
&J^{\beta}(x)J^{c}(w)=\frac{k\delta_{c,\beta}}{(x-w)^2}+\sum_c\frac{f_{\beta,c,d}J^{d}(w)}{(x-w)}+(J^{\beta}J^{c})(w)+...\\
\end{split}
\end{align}

After substitution of \eqref{jk4} into  \eqref{jk3} one gets

\begin{align}
\begin{split}\label{2}
& \frac{k\delta_{a,\beta} J^{\gamma}(w)}{(z-w)^2}+
\frac{k\delta_{a,\gamma}J^{\beta}(w)}{(z-w)^2}+\sum_c\frac{1}{2\pi
i}\oint\frac{dx}{x-w}\frac{f_{a,\beta,c}J^{c}(x)J^{\gamma}(w)}{(z-x)}+\frac{f_{a,\gamma,c}J^{\beta}(x)J^{c}(w)}{(z-w)}=
\\&=
\frac{k\delta_{a,\beta} J^{\gamma}(w)}{(z-w)^2}+
\frac{k\delta_{a,\gamma}J^{\beta}(w)}{(z-w)^2}+\\&+\sum_c\frac{1}{2\pi
i}\oint\frac{dx}{x-w}\frac{f_{a,\beta,c}(\frac{k\delta_{c,\gamma}}{(x-w)^2}+\sum_c\frac{f_{c,\gamma,d}J^{d}(w)}{(x-w)}+(J^{c}J^{\gamma})(w)+...)}{(z-x)}+
\\&+\frac{if_{a,\gamma,c}(\frac{k\delta_{c,\beta}}{(x-w)^2}+\sum_c\frac{f_{\beta,c,d}J^{d}(w)}{(x-w)}+(J^{\beta}J^{c})(w)+...)}{(z-w)}
\end{split}
\end{align}

The following expansion takes place

\begin{align}\begin{split}\label{1}
&\frac{1}{z-x}=\frac{1}{(z-w)-(x-w)}=\frac{1}{z-w}\frac{1}{1-\frac{x-w}{z-w}}=\\
&=\frac{1}{z-w}+\frac{x-w}{(z-w)^2}+\frac{(x-w)^2}{(z-w)^3}+...
\end{split}\end{align}

Substitute \eqref{1} into \eqref{2}, one gets

\begin{align}
\begin{split}
& \frac{k\delta_{a,\beta} J^{\gamma}(w)}{(z-w)^2}+
\frac{k\delta_{a,\gamma}J^{\beta}(w)}{(z-w)^2}+ \sum_c\frac{1}{2\pi
i}\oint\frac{dx}{x-w}(\frac{f_{a,\beta,c}k\delta_{c,\gamma}}{(z-w)^3}+\sum_c\frac{f_{a,\beta,c}f_{c,\gamma,d}J^{d}(w)}{(z-w)^2}+\\
&+f_{a,\beta,c}(J^{c}J^{\gamma})(w)+...)+\frac{f_{a,\gamma,c}(\frac{k\delta_{c,\beta}}{(x-w)^2}+\sum_c\frac{f_{\beta,c,d}J^{d}(w)}{(x-w)}+(J^{\beta}J^{c})(w)+...)}{(z-w)}=\\
%\sum_c \frac{1}{2\pi
%i}\oint dx\frac{1}{(z-w)^3} \frac{if_{a,\beta,c} k\delta_{c,\gamma}}{(x-w)}+\frac{1}{(z-w)^3} \frac{if_{a,\beta,c} k\delta_{c,\gamma}}{(z-x)}+\\&+\frac{1}{(z-w)^2} \frac{if_{a,\beta,c}if_{c,\gamma,d}J^d(w)}{x-w}+\frac{1}{(z-w)^2} \frac{if_{a,\beta,c}if_{c,\gamma,d}J^d(w)}{z-x}
%+\frac{1}{(z-w)} \frac{if_{a,\beta,c}(J^cJ^{\gamma})(w)}{x-w}+\\&+\frac{1}{(z-w)} \frac{if_{a,\beta,c}(J^cJ^{\gamma})(w)}{z-x}+\frac{1}{(z-w)} \frac{if_{a,\gamma,c}(J^{\beta}J^{c}(w))}{x-w}=\\&=
%\frac{k\delta_{a,\beta}
%J^{\gamma}(w)}{(z-w)^2}+
%\frac{k\delta_{a,\gamma}J^{\beta}(w)}{(z-w)^2}+\frac{if_{a,\beta,\gamma} k}{(z-w)^3}-\frac{if_{a,\beta,\gamma} k}{(z-w)^3}+\\&+\frac{if_{a,\beta,c}if_{c,\gamma,d}J^d(w)}{(z-w)^2}-\frac{if_{a,\beta,c}if_{c,\gamma,d}J^d(w)}{(z-w)^2}
%+\frac{if_{a,\beta,c}(J^cJ^{\gamma})(w)}{z-w}-\\&-\frac{if_{a,\beta,c}(J^cJ^{\gamma})(w)}{(z-w)^2}+\frac{if_{a,\gamma,c}(J^{\beta}J^{c})(w)}{z-w}=\\
%&=\frac{k\delta_{a,\gamma}J^{\beta}(w)}{(z-w)^2}++\frac{if_{a,\gamma,c}(J^{\beta}J^{c})(w)}{z-w}
&=\frac{k\delta_{a,\beta} J^{\gamma}(w)}{(z-w)^2}+
\frac{k\delta_{a,\gamma}J^{\beta}(w)}{(z-w)^2}+
\frac{f_{a,\beta,\gamma}k}{(z-w)^3}+\\&+\frac{f_{a,\beta,c}f_{c,\gamma,d}J^d(w)}{(z-w)^2}+
\frac{f_{a,\beta,c}(J^cJ^{\gamma})(w)}{z-w}+\frac{f_{a,\gamma,c}(J^{\beta}J^c)(w)}{z-w}.
\end{split}
\end{align}

\label{ap1}

\subsection{Proof of Lemma \ref{l2}} \label{ap4}

One has
%{\bf 2. The case $\frac{PfF_{F_{I_2}(J\setminus
%I_1)}(z)}{(z-x)^3}$.}

\begin{align}
\begin{split}\label{df}
& \sum_{I}\sum_{I=I_1\sqcup I_2}\frac{1}{2\pi i}\oint \frac{dx}{x-w}
\frac{PfF_{F_{I_2}(J\setminus I_1)}(z)PfF_I(w)}{(z-x)^3}=
\sum_{I}\sum_{I=I_1\sqcup I_2}\frac{PfF_{F_{I_2}(J\setminus
I_1)}(z)PfF_I(w)}{(z-w)^3}.
 \end{split}
\end{align}

Let us  consider the expression

\begin{equation}
 \sum_{I}\sum_{I=I_1\sqcup I_2}\frac{PfF_{F_{I_2}(J\setminus I_1)}(z)PfF_I(w)}{(z-w)^3}.
\end{equation}

Using the Proposition \ref{p4} we obtain that the considered product
equals

\begin{align}
\begin{split}\label{ddd}
& \sum_{I}\sum_{I=I_1\sqcup I_2}(\frac{PfF_{I\setminus
(F_{I_2}(J\setminus
I_1))}(w)}{(z-w)^5}+\frac{PfF_{F_{F_{I_2}(J\setminus
I_1)}I}(w)}{(z-w)^4}+\frac{(F_{F_{I_2}(J\setminus
I_1)}PfF_I)(w)}{(z-w)^3}+l.o.t)
 \end{split}
\end{align}

Let us prove that the term of the highest order in \eqref{ddd} is
equal to zero. Prove the following Proposition.

\begin{prop}
\label{p6} Let $|I|=4$ and let us be given a partition $I=I_1\sqcup
I_2$, $|I_1|=|I_2|=2$. Let us also be given another set of indices
$J$, such that $|J|=4$. Then
$$F_{I_2}(J\setminus I_1)$$  does not equal to zero and does not contain
repetitions if and only if $I_1\subset J$, one of the elements of
$I_2$ is contained in $J$, and the other is not contained in $J$. In
particular  $I$ and $J$ differ only by one element.
\end{prop}

\proof

Indeed if $I_1\nsubseteq J$ then $J\setminus I_1=0$. If $I_1\subset
J$ then two cases are possible.  In the first case $I_2\subset
J\setminus I_1$ (that is actually $I=J$). Then $F_{I_2}(J\setminus
I_1)$ contains repetitions.  In the second case $I_2\nsubseteq
J\setminus I_1$, then $J\setminus I_1$ contains an  one element that
does not belong to $I$. Then $F_{I_2}(J\setminus I_1)$ also contains
one  such an element.
%, that is
%$F_{I_2}(J\setminus I_1)\nsubseteq J$.
\endproof

As a corollary one get that either
 $F_{I_2}(J\setminus I_1)\nsubseteq I$ or $F_{I_2}(J\setminus I_1)$ equals to zero, or contains repetitions. In all cases  $PfF_{I\setminus (F_{I_2}(J\setminus
I_1))}(w)=0$.

Thus the terms of the highest order in \eqref{ddd} equals to zero.

%{ \bf 3. Случай $\frac{PfF_{F_{F_{I_2}(J\setminus
%I_1)}I}(w)}{(z-w)^4}$}

Consider the next term in \eqref{ddd}

\begin{align}
\begin{split}\label{df1}
& \sum_{I}\sum_{I=I_1\sqcup I_2}\frac{PfF_{F_{F_{I_2}(J\setminus
I_1)}I}(w)}{(z-w)^4}
\end{split}
\end{align}

According to Proposition \ref{p6},  $F_{I_2}(J\setminus I_1)$ is not
equal to zero and does not contain repetitions if and only if the
sets $I$ and $J$  differ by exactly one element, that is the
following equalities take place

\begin{align}
\begin{split}
J=\{a,b,c,e\},\,\,\, I=\{a,b,c,d\}\text{ or } \{a,b,d,e\} \text{ or
} \{a,d,c,e\} \text{ or } \{d,b,c,e\}.
\end{split}
\end{align}

Consider possible choices for partitions $I=I_1\sqcup I_2$ и and
calculate $F_{I_2}(J\subset I_1)$.

If  $I=\{a,b,c,d\}$ then

\begin{align}
\begin{split}
& I_1=\{a,b\},\,\,\,     I_2=\{c,d\},\,\,\,    (J\setminus I_1)=\{c,e\},\,\,\,  F_{I_2}(J\setminus I_1)=-\{d,e\}  \\
& I_1=-\{a,c\},\,\,\,     I_2=-\{b,d\},\,\,\,    (J\setminus I_1)=-\{b,e\},\,\,\,  F_{I_2}(J\setminus I_1)=-\{d,e\}  \\
& I_1=\{b,c\},\,\,\,     I_2=\{a,d\},\,\,\,    (J\setminus I_1)=\{a,e\},\,\,\,  F_{I_2}(J\setminus I_1)=-\{d,e\}  \\
\end{split}
\end{align}

If  $I= \{a,b,d,e\}$ then

\begin{align}
\begin{split}
& I_1=\{a,b\},\,\,\,     I_2=\{d,e\},\,\,\,    (J\setminus I_1)=\{c,e\},\,\,\,  F_{I_2}(J\setminus I_1)=\{c,d\}  \\
& I_1=-\{b,e\},\,\,\,     I_2=-\{a,d\},\,\,\,    (J\setminus I_1)=-\{a,c\},\,\,\,  F_{I_2}(J\setminus I_1)=-\{d,c\}  \\
& I_1=\{a,e\},\,\,\,     I_2=\{b,d\},\,\,\,    (J\setminus I_1)=\{b,c\},\,\,\,  F_{I_2}(J\setminus I_1)=-\{d,c\}  \\
\end{split}
\end{align}

If  $I=\{a,d,c,e\} $ then

\begin{align}
\begin{split}
& I_1=\{c,e\},\,\,\,     I_2=\{a,d\},\,\,\,    (J\setminus I_1)=\{a,b\},\,\,\,  F_{I_2}(J\setminus I_1)=-\{d,b\}  \\
& I_1=-\{a,c\},\,\,\,     I_2=-\{d,e\},\,\,\,    (J\setminus I_1)=-\{b,e\},\,\,\,  F_{I_2}(J\setminus I_1)=\{b,d\}  \\
& I_1=\{a,e\},\,\,\,     I_2=\{d,c\},\,\,\,    (J\setminus I_1)=\{b,c\},\,\,\,  F_{I_2}(J\setminus I_1)=\{b,d\}  \\
\end{split}
\end{align}

If  $I= \{d,b,c,e\}$ then

\begin{align}
\begin{split}
& I_1=\{c,e\},\,\,\,     I_2=\{d,b\},\,\,\,    (J\setminus I_1)=\{a,b\},\,\,\,  F_{I_2}(J\setminus I_1)=\{a,d\}  \\
& I_1=-\{b,e\},\,\,\,     I_2=-\{d,c\},\,\,\,    (J\setminus I_1)=-\{a,c\},\,\,\,  F_{I_2}(J\setminus I_1)=\{a,d\}  \\
& I_1=\{b,c\},\,\,\,     I_2=\{d,e\},\,\,\,    (J\setminus I_1)=\{a,e\},\,\,\,  F_{I_2}(J\setminus I_1)=\{a,d\}  \\
\end{split}
\end{align}

%В первом случае имеем

%\begin{align}
%\begin{split}
%& I_1 ,\,\,\,I_2,\,\,\,J\setminus I_1,\,\,\, F_{I_2}(J\setminus I_1)\\
%&\{a,b\},\{c,d\},\{c,e\},-\{d,e\}\\
%&-\{a,c\},-\{b,d\},-\{b,e\},-\{d,e\}\\
%&\{b,c\},\{a,d\},\{a,e\},-\{d,e\}
%\end{split}
%\end{align}

Thus in the case $I= \{a,b,c,d\}$ one has

\begin{align}
\begin{split}
\frac{PfF_{-F_{\{d,e\}}\{a,b,c,d\}}(w)}{(z-w)^4}=\frac{PfF_{\{a,b,c,e\}}(w)}{(z-w)^4}=\frac{PfF_{J}(w)}{(z-w)^4}
\end{split}
\end{align}

In the other cases considerations are similar. Thus  \eqref{df1}
equals to

\begin{equation}
\frac{12(N-4)PfF_{J}(w)}{(z-w)^4}.
\end{equation}

Here the factor $12$ in the numerator is a product of the factor $3$
that corresponds  corresponds to the choice of one of four case
described above  and the factor $4$ that corresponds to the choice
of the partition $I=I_1\sqcup I_2$. The factor $(N-4)$ corresponds
to the choice an element $d\in I$, that does not belong to $J$.

Consider the third term in  \eqref{ddd}

\begin{align}
\begin{split}
\label{dwer}
 &\sum_{I}\sum_{I=I_1\sqcup
I_2}\frac{(PfF_{F_{I_2}(J\setminus I_1)}PfF_I)(w)}{(z-w)^3}.
\end{split}
\end{align}

 According to Proposition \ref{p6} a summand in this expression can be nonzero only in the for cases that were
 described above. After summation of inputs from these cases one
 gets

\begin{align}
\begin{split}\label{pfk}
&\frac{1}{(z-w)^4}(-(F_{de}PfF_{\{a,b,c,d\}})(w)+(F_{cd}PfF_{\{a,b,d,e\}})(w)-\\&-(F_{db}PfF_{\{a,d,c,e\}})(w)+(F_{ad}PfF_{\{d,b,c,e\}})(w_).
\end{split}
\end{align}

The numerator in this expression can be simplified. Let us write the
terms of the first and the second pfaffians corresponding to
identical permutation and the permutation that changes the first and
the third, the second and the fourth indices.

\begin{align}
\begin{split}
&-(F_{de}(F_{a,b}F_{c,d}))(w)-(F_{de}(F_{c,d}F_{a,b}))(w)+(F_{c,d}(F_{a,b}F_{d,e}))(w)+(F_{c,d}(F_{d,e}F_{a,b}))(w).
\end{split}
\end{align}
Let us simplify  the sum of the first and the fourth summands:
%, the
%second and the fourth summands. For the first and the third

\begin{align}
\begin{split}
\label{summ}
&-(F_{de}(F_{a,b}F_{c,d}))(w)+(F_{c,d}(F_{d,e}F_{a,b}))(w).
\end{split}
\end{align}

Let us transform the first summand. We are going to use the
equalities \cite{Di}.

\begin{align}
\begin{split}\label{ab1}
& ((AB)E)(w)-(A(BE))(w)=(A([E,B]))(w)+(([E,A])B)(w)+([(AB),E])(w),
\end{split}
\end{align}
\begin{align}
\begin{split}\label{ab2} &
(BA)(w)=(AB)(w)-\partial
(AB)_{-1}(w)+\frac{1}{2!}\partial^2(AB)_{-2}(w)-\frac{1}{3!}\partial^3(AB)_{-3}(w)+..
\end{split}
\end{align}

 Using \eqref{ab1} one gets

\begin{align}
\begin{split}
&-(F_{de}(F_{a,b}F_{c,d}))(w)=-((F_{de}F_{a,b})F_{c,d})(w),
\end{split}
\end{align}

using \eqref{ab2}, one gets

\begin{align}
\begin{split}
&-((F_{de}F_{a,b})F_{c,d})(w)=-(F_{c,d}(F_{de}F_{a,b}))(w)+\partial(F_{c,d}(F_{de}F_{a,b}))_{-1}(w)-\\&
-\frac{1}{2!}\partial^2(F_{c,d}(F_{de}F_{a,b}))_{-2}(w)+...
\end{split}
\end{align}

Thus \eqref{summ} equals to

\begin{align}
\begin{split}
&\partial(F_{c,d}(F_{de}F_{a,b}))_{-1}(w)-\frac{1}{2!}\partial^2(F_{c,d}(F_{de}F_{a,b}))_{-2}(w)+...
\end{split}
\end{align}

By symmetry it also equals

\begin{align}
\begin{split}
&\partial(F_{d,e}(F_{a,b}F_{c,d}))_{-1}(w)-\frac{1}{2!}\partial^2(F_{d,e}(F_{a,b}F_{c,d}))_{-2}(w)+...
\end{split}
\end{align}

 Take into account other terms of the
pfaffian in \eqref{pfk}. One gets that  \eqref{pfk} equals

\begin{align}
\begin{split}\label{pfk1}
&\frac{1}{(z-w)^4}(-\partial(F_{de}PfF_{\{a,b,c,d\}})(w)+\partial(F_{cd}PfF_{\{a,b,d,e\}})(w)-\\&-\partial(F_{db}PfF_{\{a,d,c,e\}})(w)+\partial(F_{ad}PfF_{\{d,b,c,e\}})(w_).
\end{split}
\end{align}

   Above in the proof of this Lemma it was
shown that $-(F_{d,e}PfF_{\{a,b,c,d\}})_{-k}(w)=0$ for  $k>1$ and
$-(F_{d,e}PfF_{\{a,b,c,d\}})_{-1}(w)=PfF_{\{a,b,c,d\}}(w)$.
Analogously all other pfaffians are considered. Thus the sum of the
first and the second pfaffians in \eqref{pfk} equals

\begin{align}
\begin{split}
&12(N-4)\partial PfF_{\{a,b,c,e\}}(w)=12(N-4)\partial PfF_{J}(w).
\end{split}
\end{align}

Thus the expression \eqref{dwer} equals

\begin{equation}
\frac{12(N-4)\partial PfF_J(w)}{(z-w)^3}
\end{equation}

\subsection{The proof of Lemma \ref{l3}}
\label{ap5}

%{ \bf 3. The  third case.}

It is enough to prove that  when one calculates the expansion of the
expression

\begin{align}
\begin{split}\label{65}
&\frac{(F_{I_2}PfF_{J\setminus I_1})(z)PfF_I(w)}{(z-w)^2}
\end{split}
\end{align}

 one does not obtain singular terms of orders greater than $(z-w)^{-2}$.  The expression
 \eqref{65}  is nonzero only if two element from $I$ belong to $J$,
that is these
 sets are of type

$$I=\{a,b,c,d\},\,\,\,J=\{a,b,e,f\}.$$

The numerator in \eqref{65} is of type

\begin{align}
\begin{split}
\label{66} &(F_{cd}F_{ef})(z)PfF_{\{a,b,c,d\}}(w).
\end{split}
\end{align}

Using notation $I'_1=\{c,d\}$, $I'_2=\{e,f\}$, $I=\{a,b,c,d\}$ let
us write \eqref{66} as follows

\begin{align}
\begin{split}
(F_{I'_1}F_{I'_2})(z)PfF_{I}(w).
\end{split}
\end{align}

Apriory in the calculation of this OPE the terms of the order
$(z-w)^{-4}$,  $(z-w)^{-3}$,  and terms of lower order. Let us show
that actually the terms of  orders $4$ and $3$ do not appear.

Indeed, the OPE is calculated using Proposition  \ref{p4}. But since
$J\setminus I_1\nsupseteq I_1$, then $I\neq I'_1\cup I'_2$ and the
terms of the highest order in the OPE vanishes.

The second terms corresponds to

\begin{align}
\begin{split}
& F_{I'_1}(I\setminus I'_2)=F_{cd}\{a,b\}\delta_{(ab),(cd)}=0,\,\,\,
F_{I'_2}(I\setminus I'_1)=F_{ef}\{a,b\}=0
\end{split}
\end{align}

Thus this terms also vanishes. Thus there are no terms of the order
$4$ and $3$.

\subsection{Proof of Lemma \ref{l5}}
\label{ap6}

The expression on the right hand side in \eqref{167} equals

\begin{align}
\begin{split}
\label{67} &\sum_{I=I_1\sqcup
 I_2,|I_1|=|I_2|=2}(-1)^{(I_1,I_2)}
 \frac{(F_{I_1}PfF_{F_{I_2}J})(z)PfF_I(w)}{(z-w)}.
\end{split}
\end{align}

Let us divide proof of the lemma into several steps.  We shall
investigate the OPE in the numerator in \eqref{67}.

{\bf Step 1.} Show that the OPE in the numerator begins with the
term of order $(z-w)^{-3}$.  Indeed this OPE can be calculated using
Proposition \ref{p4}. Since $|I|=4$,  than the most singular term
that one can apriory obtain is proportional to $(z-w)^{-4}$.  But a
term of order $(z-w)^{-4}$ can appear only if
$$I\subset I_1\cup F_{I_2}J\Leftrightarrow I_2\subset F_{I_2}J.$$

This inclusion is impossible.  Thus after calculation of the OPE in
the numerator in \eqref{67}, one gets an expansion that begins with
$(z-w)^{-3}$.

{\bf Step 2.} Let us find a relations between a coefficient at
$(z-w)^{-3}$ and a coefficient at $(z-w)^{-2}$ in the OPE in the
numerator in \eqref{67}. Let us show that the coefficient at
$(z-w)^{-2}$ is a derivative of the coefficient of the coefficient
at $(z-w)^{-3}$.

 For this purpose put $a(z)=([PfF_I,PfF_J])(z)$, $b(w)=PfF_I(w)$.
 Firstly show that $(ab)_{-3}(w)$ equals to a coefficient at
 $(z-w)^{-4}$, and $(ab)_{-2}(w)$ equals to a coefficient at
 $(z-w)^{-3}$. Then show that $\partial (ab)_{-3}(w)=(ab)_{-2}(w)$

%Consider the expression

%\begin{align}
%\begin{split}\label{pfa}
%\sum_{I_1,I_2}(-1)^{(I_1,I_2)}\frac{(F_{I_1}PfF_{F_{I_2}J})(z)PfF_{I}(w)}{(z-w)}+\frac{PfF_I(z)(F_{I_1}PfF_{F_{I_2}J})(w)}{(z-w)}.
%\end{split}
%\end{align}

%Show that in this expression there appear no terms of orders
%$(z-w)^{-4}$, $(z-w)^{-3}$.

Let us use the following fact.  In \cite{AG} a formula for a
commutator of two pfaffians was proved. In the case of two pfaffians
of order $4$ it looks as follows

\begin{align}
\begin{split}\label{com}
[PfF_I,PfF_J]=PfF_{PfF_IJ}+\sum_{I=I_1\sqcup
I_2}(-1)^{(I_1,I_2)}F_{I_1}PfF_{F_{I_2}J}.
\end{split}
\end{align}

Thus one has

\begin{equation}
a(z)b(w)=PfF_{PfF_IJ}(z)PfF_I(w)+\sum_{I=I_1\sqcup
I_2}(-1)^{(I_1,I_2)}(F_{I_1}PfF_{F_{I_2}J})(z)PfF_I(w).
\end{equation}

To do the first step one must prove the Proposition

\begin{prop}
The OPE
\begin{equation}
\label{45} PfF_{PfF_IJ}(z)PfF_I(w)
\end{equation}
 begins with the term $(z-w)^{-1}$.
\end{prop}

\proof To calculate the OPE in the numerator the formula
\eqref{long} is used.

In this calculation a terms of order $4$ can appear. More precise it
appears if
$$PfF_IJ=I.$$

 Note that the number of elements of $J$, that contain in $I$, is the same as the number of elements of $PfF_IJ$,
 that contain in $I$. Thus the equality written above is possible if and only if $J=I$. But then $PfF_IJ$
 contains repetitions. Thus the terms of order $4$  does not appear.

The terms of order $3$ appears if
$$F_{I_2}((PfF_IJ)\setminus I_1)$$ is nonzero and does not contain repetitions for some partition  $I=I_1\sqcup
I_2$. Put $\mathcal{J}=PfF_IJ$. Above it was shown that
$F_{I_2}(\mathcal{J}\setminus I_1)$ has this property if and only if
$\mathcal{J}$ and $I$ differ by exactly one element. If $PfF_IJ$ and
$I$ differ by one element then $J$ and $I$ also differ by one
element. In this case $PfF_IJ$ contains repetitions. Thus the term
of the order $3$ does not appear.

The term of the order $2$ appears if
$$(F_{I_2}PfF_{PfF_IJ\setminus I_1})(w)\neq 0$$  for some partition $I=I_1\sqcup I_2$.
As above one can prove that this is impossible. Thus the terms of
order $2$ do not appear.

Thus the expansion of the expression \eqref{45} begins with the term
$(z-w)^{-1}$.

\endproof

Let us show that $\partial (ab)_{-3}(w)=(ab)_{-2}(w)$. In \cite{AG1}
it was shown that the tensor

\begin{equation}\label{ei}
\sum_{\{i_1,...,i_4\}}(e_{i_1}\wedge...\wedge e_{i_4})\otimes
(e_{i_1}\wedge...\wedge e_{i_4}))
\end{equation}

is invariant under the action of $\mathfrak{o}_N$. Here
$i_k=1,...,N$ and $e_i$ is a standard base in $\mathbb{C}^N$.

From the invariance of \eqref{ei} it follows that

\begin{align}
\begin{split}\label{pfb}
([PfF_I,PfF_J])(z)PfF_I(w)+PfF_I(z)([PfF_I,PfF_J])(w)=0.
\end{split}
\end{align}

Thus one has

\begin{equation}
\label{opesh}
 a(z)b(w)+b(z)a(w)=0.
\end{equation}

Take OPEs in \eqref{opesh}, one gets that

\begin{equation}
(ab)_{-2}(w)=\partial(ab)_{-3}(w).
\end{equation}

Thus the coefficient at $(z-w)^{-2}$ equals to the coefficient at
$(z-w)^{-3}$.

{\bf Step 3.} Let us calculate a coefficient at $(z-w)^{-3}$ in the
expansion of the numerator in \eqref{67}, that is a coefficient at
$(z-w)^{-3}$ in the expansion of

\begin{equation}
\label{raz} \sum_{I}\sum_{I=I_1\sqcup
I_2}(F_{I_1}PfF_{F_{I_2}J})(z)PfF_I(w).
\end{equation}

To calculate  this coefficient let us write

\begin{align}
\begin{split}
\label{i1i2} &\sum_{I=I_1\sqcup
I_2}(-1)^{(I_1,I_2)}(F_{I_1}\overbrace{PfF_{F_{I_2}J})(z)PfF_I(w)}=\sum_{I=I_1\sqcup
I_2}(-1)^{(I_1,I_2)}PfF_I(w)(F_{I_1}PfF_{F_{I_2}J})(z)=\\
&=\sum_{I=I_1\sqcup I_2}(-1)^{(I_1,I_2)}\frac{1}{2\pi i} \oint
\frac{dx}{x-w}
\overbrace{PfF_{I}(w)F_{I_1}(x)}PfF_{F_{I_2}}(z)+F_{I_1}(x)\overbrace{PfF_{I}(w)PfF_{F_{I_2}}(z)}=\\
&=\sum_{I_1,I_2}(-1)^{(I_1,I_2)}\oint \frac{dx}{x-w}
\overbrace{F_{I_1}(x)PfF_{I}(w)}PfF_{F_{I_2}}(z)+F_{I_1}(x)\overbrace{PfF_{F_{I_2}}(z)PfF_{I}(w)}.
\end{split}
\end{align}

Consider two summands in \eqref{i1i2} separately. Let us substitute
the expression for the contraction
$\overbrace{F_{I_1}(x)PfF_{I}(w)}$ which is given by Proposition
\ref{p4}. After an explicit calculation one gets that the term of
order $(z-w)^{-3}$  equals to

\begin{equation}
\sum_{I_1,I_2}(-1)^{(I_1,I_2)}\frac{(k+2)PfF_{PfF_{I\setminus
I_1}F_{I_2}J}}{(z-w)^3}.
\end{equation}

Note that $PfF_{I\setminus I_1}=(-1)^{(I_1,I_2)}F_{I_2}$. Thus the
numerator  equals $-PfF_J$, if exactly one of the elements of $I_2$
belongs to $J$ and equals to zero otherwise. There are

\begin{equation}
\frac{4(N-4)}{2}\frac{(N-2)(N-3)}{2}
\end{equation}
pair of such subsets $I_1,I_2\subset\{1,...,N\}$.

Indeed the factor $4$ corresponds to the choice of the element of
$I_2$, that does not belong to $J$, the factor $(N-4)$ corresponds
to the choice of the second element, the denominator $2!$
corresponds to the fact that the set $I_2$ is unordered. The factor
$\frac{(N-2)(N-3)}{2!}$ is the number of choices of the set $I_1$.

Thus the first summand in \eqref{i1i2} equals

\begin{equation}
\label{pervo} -\frac{(N-2)(N-3)(N-4)(k+2)PfF_J(w)}{(z-w)^3}.
\end{equation}

Consider the second summand in \eqref{i1i2}. The contraction
$\overbrace{PfF_{F_{I_2}}(z)PfF_{I}(w)}$ is given by Proposition
\ref{p5}. The numerators in this expression depend on $z$. After
integration $\oint \frac{dx}{x-w}$  one gets an expression where the
term of order $(z-w)^{-3}$ is given by expression

\begin{equation}\label{vyr}
\sum_{I_1,I_2}(-1)^{(I_1,I_2)}\sum_{I'_1,I'_2}(-1)^{(I'_1,I'_2)}\frac{-2(k+2)(F_{I_1}PfF_{F_{I'_2}((F_{I_2}J)\setminus
I'_1)})(w)}{(z-w)^3}.
\end{equation}

Consider firstly the expression

\begin{equation}
\label{su}
\sum_{I'_1,I'_2}(-1)^{(I'_1,I'_2)}(F_{I_1}PfF_{F_{I'_2}((F_{I_2}J)\setminus
I'_1)})(w)
\end{equation}

Let us do the summation over partitions $I=I'_1\sqcup I'_2$. By
Proposition \ref{p6}, $F_{I'_2}((F_{I_2}J)\setminus I'_1)\neq 0$ if
and only if $I'_1\subset F_{I_2}J$, and one of the elements $I'_2$
belongs to $F_{I_2}J$, and the other does not belong.
% For z fixed
%set $I$  the number of such partitions equals to $3$.

For each such partition, doing the same calculation as in the proof
of Lemma \ref{l2} (using everywhere $\mathcal{J}=F_{I_2}J$ instead
of $J$), one gets that

\begin{equation}
\label{wer} \sum_{I=I'_1\sqcup
I'_2}(F_{I_1}PfF_{F_{I'_2}((F_{I_2}J)\setminus
I'_1)})(w)=3(F_{I_1}F_{I\Delta (F_{I_2}J)})(w),\footnote{The set $I$
is fixed}
\end{equation}

where $I\Delta (F_{I_2}J)$ is a symmetric difference of sets. Thus
for $I=\{a,b,c,d\}$ and $F_{I_2}J=\{a,b,c,e\}$ one has $I\Delta
(F_{I_2}J)=-\{d,e\}$.

 Find the sum $\sum_I\sum_{I=I_1\sqcup
I_2}$ of the expressions  \eqref{wer}.  Let $J=\{a,b,c,e\}$. The set
$I_1$ must be a subset in $J$, one of the elements of $I_2$ must be
contained in $J$ and the other not.  Thus the possibilities for
$I_1$ and $I_2$ are the following

\begin{align}
\begin{split}
& I_1=\{a,b\},\,\,\,     I_2=\{с,d\},\,\,\,      F_{I_2}J=-\{a,b,d,e\},\,\,\,   I\Delta (F_{I_2}J)=-\{c,e\} \\
& I_1=-\{a,c\},\,\,\,     I_2=-\{b,d\},\,\,\,    F_{I_2}J=-\{a,d,c,e\},\,\,\,   I\Delta (F_{I_2}J)=-\{b,e\}\\
& I_1=\{b,c\},\,\,\,     I_2=\{a,d\},\,\,\,   F_{I_2}J=-\{d,b,c,e\},\,\,\,    I\Delta (F_{I_2}J)-\{a,e\}\\
\end{split}
\end{align}

Thus the sum of expressions \eqref{wer} equals

\begin{equation}
-3PfF_J(w)
\end{equation}

Hence the second summand in \eqref{i1i2} equals

\begin{equation}
\frac{-6(k+2)PfF_J(w)}{(z-w)^3}
\end{equation}

So, the expression \eqref{i1i2} equals

\begin{equation}
-\frac{(k+2)((N-2)(N-3)(N-4)+6)PfF_J(w)}{(z-w)^3}.
\end{equation}

%By Proposition \ref{p6}, $F_{I'_2}((F_{I_2}J)\setminus I'_1)\neq 0$
%if and only if $F_{I_2}J$  and $I$ differ by only one element. This
%is possible if and only if $I_1\subset J$, and one of the elements
%of $I_2$ belongs to $J$, and another does not belong. Doing the same
%calculation as in the proof of Lemma \ref{l2}, one gets that

%\begin{equation}
%\sum_{I'_1,I'_2}(-1)^{(I'_1,I'_2)}(F_{I_1}PfF_{F_{I'_2}((F_{I_2}J)\setminus
%I'_1)})(w)=?(F_{I_1}F_{J\setminus I_2})(w).
%\end{equation}

%Thus after summation $\sum_{I_1,I_2}(-1)^{(I_1,I_2)}$ one gets that
%the expression \eqref{su} equals

%\begin{equation}
%\frac{-2(k+2)PfF_J(w)}{(z-w)^3}+l.o.t.
%\end{equation}

{\bf Conclusion.} Thus
%the expansion of the expression \eqref{raz}
%is the following

%\begin{align}
%\begin{split}
%\frac{(k+2)((N-2)(N-3)(N-4)-6)PfF_J(w)}{(z-w)^3}+l.o.t.
%\end{split}
%\end{align}

%Hence
the coefficient at $(z-w)^{-4}$ in the expansion of \eqref{67}
equals
$$-(k+2)((N-2)(N-3)(N-4)+6)PfF_J(w),$$

the coefficient at $(z-w)^{-3}$ in the expansion of \eqref{67}
equals
$$-(k+2)((N-2)(N-3)(N-4)+6)\partial PfF_J(w).$$

\end{document}